%% file: main.tex
\documentclass[12pt]{article}

\include{prematter}

\title{\Large Certified Lumped Approximations \\for the Conduction Dunking Problem }

\author[1]{\bf{Kento Kaneko}}
\author[2]{\bf{Claude Le Bris}}
\author[1]{\bf{Anthony T Patera}}

\affil[1]{\small{Massachusetts Institute of Technology, Cambridge, MA 02139, USA. \newline\texttt{kaneko@mit.edu}, \texttt{patera@mit.edu}}}

\affil[2]{\'Ecole des Ponts and INRIA, Champs-sur-Marne, 77455 Marne La Vall\'ee, France. \newline\texttt{claude.le-bris@enpc.fr}}

\date{}

\begin{document}

\maketitle
\vspace{-0.9cm}
\centerline{\bf \small Acknowledgments}
\vspace{0.15cm}

{\small
\small We thank Prof. Masayuki Yano for his many contributions to our finite element software in particular related to adaptive refinement and associated {\it a posteriori} error estimation.

This work is supported by ONR Grant N00014-21-1-2382, Grant Monitor Dr Reza Malek-Madani. The research of the second author is partially supported by ONR and EOARD; part of the work was completed while the second author was visiting MIT.
}

\input{abs}      
\input{intro}    
\input{form}     
\input{fcnals}   
\input{result}   
\input{numer}    
\input{lump}     
\input{phi}      

{
  \bibliographystyle{ieeetr}
  \addcontentsline{toc}{part}{Bibliography}
  \bibliography{sbp}
}

\end{document}

%% file: prematter.tex
\usepackage{url}
\usepackage{breakurl}
\usepackage[breaklinks]{hyperref}

\usepackage[left=1.5in,right=1.5in,top=1.in,bottom=1.in]{geometry}
\usepackage{tikz}
\usepackage{authblk}
\usetikzlibrary{arrows.meta, positioning, calc}
\usepackage{verbatim,setspace}
\usepackage{lmodern}
\usepackage{pdfcomment}

\usepackage{setspace}
\usetikzlibrary{circuits.ee.IEC, positioning}

\usepackage[rflt]{floatflt}
\usepackage{float}
\usepackage{subcaption}

\usepackage{graphicx}
\graphicspath{{./Figures/}}

\usepackage{amsthm,amsmath,amssymb,amsfonts,amscd,amsbsy,thmtools}
\usepackage{epstopdf}
\usepackage{calrsfs}
\DeclareMathAlphabet{\pazocal}{OMS}{zplm}{m}{n}

\usepackage{textcomp}

\usepackage{caption}
\DeclareMathAlphabet{\pazocal}{OMS}{zplm}{m}{n}

\def\OmMeas{{|\Omega |}}

\def\tf{{t_{\text{final}}}}

\def\R0{{\text{RHE}^0}}

\def\calL{\mathcal{L}}

\def\calR{\mathcal{R}}

\def\DD{\mathbb{D}}

\def\FF{\mathbb{F}}

\def\RR{\mathbb{R}}

\def\pcalO{\pazocal{O}}

\def\pcalD{\pazocal{D}}

\def\avg{{\text{avg}}}
\def\paavg{{\partial{\text{avg}}}}

\def\tf{{t_{\text{final}}}}

\def\Xint#1{\mathchoice
{\XXint\displaystyle\textstyle{#1}}
{\XXint\textstyle\scriptstyle{#1}}
{\XXint\scriptstyle\scriptscriptstyle{#1}}
{\XXint\scriptscriptstyle\scriptscriptstyle{#1}}
\!\int}
\def\XXint#1#2#3{{\setbox0=\hbox{$#1{#2#3}{\int}$ }
\vcenter{\hbox{$#2#3$ }}\kern-.6\wd0}}

\def\dashint{\Xint-}

\newcommand{\dv}[1]{{\underline{#1}}}

\def\dvt{{\dv{t}}}

\def\dvtf{{\dv{t}_\text{final}}}
\def\dvT{{\dv{T}}}
\def\dvTi{{\dv{T}_\text{i}}}
\def\dvTinf{{\dv{T}_{\infty}}}
\def\dvOm{{\dv{\Omega}}}
\def\dvpOm{{\dv{\partial\Omega}}}

\def\bnabla{\boldsymbol{\nabla}}

\usepackage{mdframed}

\definecolor{light-gray}{gray}{0.90}
\definecolor{light-blue}{RGB}{173,216,230}

\declaretheoremstyle[
spaceabove=2pt,
spacebelow=9pt,
notefont=\normalfont,
notebraces={}{},
headfont=\itshape,
bodyfont=\normalfont,
headformat=\NAME~{\em {\NUMBER}}{\normalfont }\NOTE,
mdframed={backgroundcolor=light-gray}
]{mystyle0}

\declaretheoremstyle[
spaceabove=2pt,
spacebelow=9pt,
notefont=\normalfont,
notebraces={}{},
headfont=\itshape,
bodyfont=\normalfont,
headformat=\NAME~{\em {\NUMBER}}{\normalfont }\NOTE,
mdframed={backgroundcolor=white}
]{mystyle1}

\declaretheoremstyle[
spaceabove=2pt,
spacebelow=9pt,
notefont=\normalfont,
notebraces={}{},
headfont=\itshape,
bodyfont=\normalfont,
headformat=\NAME~{\em {\NUMBER}}{\normalfont }\NOTE,
mdframed={backgroundcolor=light-blue}
]{mystyle3}

\declaretheoremstyle[
spaceabove=2pt,
spacebelow=9pt,
notefont=\normalfont,
notebraces={}{},
headfont=\itshape,
bodyfont=\normalfont,
headformat=\NAME~{\em {\NUMBER}}{\normalfont }\NOTE,
mdframed={backgroundcolor=light-gray}
]{mystyle}

\declaretheoremstyle[
spaceabove=2pt,
spacebelow=9pt,
notefont=\normalfont,
notebraces={}{},
headfont=\itshape,
bodyfont=\normalfont,
headformat=\NAME~{\em {\NUMBER}}{\normalfont }\NOTE,
mdframed={backgroundcolor=white}
]{mystyle}

\def\OmMeas{|\Omega|}

\def\psipz{\psi'^{\,0}}

\def\Rzp{\RR_{+}}

\newcommand{\Bdunk}{\text{Bi}_\text{dunk}}
\newcommand{\Bdunkp}{\text{Bi}'_\text{dunk}}

\newcommand{\uavgL}[1]{\tilde{u}_\avg^{#1}}

\newcommand{\uavgLP}{\tilde{u}_\avg^{2\text{P}}}

\newcommand{\uDeltaL}{\tilde{u}_\Delta^{2\text{P}}}

\newcommand{\taueq}{\tilde{\tau}^1_\avg}

\newcommand{\eavgL}[1]{e_\avg^{#1}}
\newcommand{\eavgLa}[1]{e_\avg^{{#1}\,\text{asymp}}}
\newcommand{\eavgLP}{e_\avg^{2\text{P}}}
\newcommand{\eavgLUB}[1]{e_\avg^{{#1}\,\text{UB}}}
\newcommand{\eDeltaLr}{e_\Delta^{{2}\,\text{rel}}}

\newcommand{\Lcond}{\calL_\text{cond}}

\definecolor{kk}{RGB}{191, 0, 255}

\DeclareMathOperator*{\essinf}{ess\,inf}

\newcommand{\sci}[2]{#1\;\!\cdot\;\!10^{#2}}

%% file: abs.tex
\begin{abstract}
We consider the dunking problem: a solid body at uniform temperature $T_\text{i}$ is placed in a environment characterized by farfield temperature $T_\infty$ and time-independent spatially uniform heat transfer coefficient; we permit heterogeneous material composition. The problem is described by a heat equation with Robin boundary conditions. The crucial parameter is the Biot number, a nondimensional heat transfer coefficient; we consider the limit of small Biot number.

We introduce first-order and second-order asymptotic approximations (in Biot number) for the spatial domain average temperature as a function of time; the first-order approximation is the standard `lumped model'. We provide asymptotic error estimates for the first-order and second-order approximations for small Biot number, and also, for the first-order approximation, non-asymptotic bounds valid for all Biot number. We also develop a second-order approximation and associated asymptotic error estimate for the normalized difference in the domain average and boundary average temperatures. Companion numerical solutions of the heat equation confirm the effectiveness of the error estimates for small Biot number.

The second-order approximation and the first-order and second-order error estimates depend on several functional outputs associated with an elliptic partial differential equation; the latter can be derived from Biot-sensitivity analysis of the heat equation eigenproblem in the limit of small Biot number. Most important is the functional output $\phi$, the only functional output required for the first-order error estimate and also the second-order approximation; $\phi$ admits a simple physical interpretation in terms of conduction length scale. We characterize a class of spatial domains for which the standard lumped-model criterion --- Biot number (based on volume-to-area length scale) small --- is deficient.
\end{abstract}


%% file: intro.tex
\section{Introduction}\label{sec:intro}

The dunking problem is ubiquitous in heat transfer engineering education and professional practice: a body characterized by spatial domain $\dvOm$ and boundary $\dvpOm$, at initial uniform temperature $\dvTi$, is immersed at time $\dvt=0$ in an environment (fluid and enclosure) at farfield temperature $\dvTinf$ \cite{LandL, Cengel, IandD}. The temperature distribution of the body, $\dvT$, evolves over the time interval $0 < \dvt\le \dvtf$. 

In this introduction, for simplicity of exposition, we largely consider homogeneous body composition, and hence uniform thermophysical properties --- in particular volumetric specific heat, $\dv{\rho c}$, and scalar thermal conductivity, $\dv{k}$; however, in subsequent development, we shall treat the important case of material heterogeneity (and hence non-uniform thermophysical properties), particularly in the form of several isotropic materials. The thermal environment is characterized, most simply, by a uniform heat transfer coefficient, $\dv{h}$, independent of time and space. 

We associate to the body domain two length scales: a chosen extrinsic length scale, $\dv{\ell}$; an intrinsic length scale, $\dv{\calL} \equiv |\dvOm|/|\dvpOm|$, where $|\dvOm|$ and $|\dvpOm|$ denote the volume and surface area of the body, respectively. We may then introduce the nondimensional variables which shall inform this work. Lengths are scaled according to the problem-defined length scale, $\dv{\ell}$, and hence our nondimensional spatial coordinate is given by $x\equiv \dv{x} / \dv{\ell}$. Times are scaled with the diffusive time scale,  ${\dv{\rho c}\,\dv{\ell}^2}/{\dv{k}}$, and hence our nondimensional temporal coordinate (Fourier number) is given by $t \equiv {\dvt\,\dv{k}}/{(\dv{\rho c}\dv{\ell}^2)}$. Temperature is nondimensionalized in standard fashion with respect to the initial temperature and farfield temperature: $u\equiv \frac{\dvT-\dvTinf}{\dvTi-\dvTinf}$. Finally, the Robin coefficient is the Biot number, $B\equiv {\dv{h}\,\dv{\ell}}/{\dv{k}}$; we also introduce the Biot number based on the intrinsic length scale, $\Bdunk = B\, \frac{\dv{\calL}}{\dv{\ell}}$. Two notational remarks: ``$\equiv$'' in this work signifies ``is defined as'' (equivalent to ``:='' ); dimensional quantities shall be indicated with underline  (and nondimensional quantities will be unadorned).

We study in this paper the dunking problem in the limit of small Biot number. The small-Biot limit arises quite often in practice, in particular for ``everyday'' artifacts at modest temperatures subject to natural convection and radiation heat transfer. In contrast, larger systems at very elevated temperatures subject to brisk forced convection or change-of-phase heat transfer typically yield larger Biot number. An example of a small-Biot application is the process of annealing by natural convection --- for instance, to relieve residual stresses.

We consider, in this small-Biot limit, three Quantities of Interest, or QoI. The first  QoI is the domain average temperature as a function of time, $\dvT_\avg (\dvt)$ --- also directly related to the heat loss (or gain) of the body from (or to) the environment. (In actual fact, we weigh the average by the volumetric specific heat normalized by the average volumetric specific heat; we present the precise definition, in the nondimensional context, in Section \ref{sec:form}.)  The second QoI is the boundary average temperature as a function of time, $\dvT_\paavg(\dvt)$.  The associated nondimensional counterparts of these first two QoI are denoted $u_\avg(t)$ and $u_\paavg(t)$, respectively. The third QoI is the normalized domain average-boundary average temperature difference as a function of time: in nondimensional form, $u_\Delta(t)\equiv\frac{u_\avg(t)-u_\paavg(t)}{u_\avg(t)}$; this QoI is a measure of the relative spatial temperature variation within the body.

The first-order\footnote{Note in general ``order'' shall refer to the convergence rate in $B$ (or, equivalently, $\Bdunk$): linear in $B$ for first-order, and quadratic in $B$ for second-order.} `classical' lumped approximation to $u_\avg$ for the small-$\Bdunk$ dunking problem is a simple exponential in time, $\uavgL{1}(t) \equiv \exp(-B \gamma t)$ \cite{LandL,Cengel,IandD}, where $\gamma \equiv 1/\calL$ and $\calL = \dv{\calL}/\dv{\ell}$. (In the case of heterogeneous material composition, $\dv{\rho c}$ is replaced by the spatial average of $\dv{\rho c}$ over $\dv{\Omega}$.) However, and despite the simplicity and utility of this small-$\Bdunk$ lumped approximation, there is relatively little analysis of the associated approximation error. For example, most textbooks \cite{LandL,Cengel,IandD} indicate only that $\Bdunk$ must be small; in some cases, a threshold might be provided, for example $\Bdunk \leq 0.1$ — but typically qualified by ``usually''. A notable exception is the important work of Gockenbach and Schmidtke \cite{ODD}, in which a rigorous asymptotic error estimate is provided for the particular case of a sphere with homogeneous thermophysical properties: $|\uavgL{1}(t) - u_\avg(t)| \leq (3/5)\,\Bdunk/e + \pcalO(B^2)$ for all time $t$. Gockenbach and Schmidtke also propose a second-order lumped approximation, again for the homogeneous sphere. The analysis of \cite{ODD} is based on an eigenvalue Taylor-series expansion in Biot.

In this work, we emphasize the mathematical results derived from rigorous analysis of the underlying parabolic partial differential equation (PDE) -- the Heat Equation -- associated with the dunking problem. Numerical results supplement the mathematical results with concrete examples. The foundation for this new analysis framework is described in \cite{KandLBandP}, which includes all necessary mathematical details as well as sketches of the proofs of all propositions. (We note that many of the mathematical results in \cite{KandLBandP} are either elementary and/or well-known --- and included only for the convenience of the reader.) We shall refer often to \cite{KandLBandP}; however, the current paper is self-contained as regards the material  important for application to engineering practice and education. We shall typically refer to \cite{KandLBandP} in expanded format as \cite{KandLBandP} [\dots], where the [\dots] will indicate the specific point in \cite{KandLBandP} at which the result is presented.

We enumerate here the specific contributions of this work:
\begin{enumerate}
\item We introduce a framework for small-$B$ approximation and error estimation formulated in terms of a perturbation field $\xi$; we can advantageously cast $\xi$ as the solution to an auxiliary elliptic PDE  derived from sensitivity analysis \cite{Joseph,Kato} of the small-$B$ eigenproblem associated to our heat equation \cite{KandLBandP} [Appendix A.2]. The key quantity in our analysis is $\phi\equiv \int_\Omega \kappa \nabla \xi \cdot\nabla \xi$, where $\kappa$ is the nondimensional thermal conductivity, defined as $\dv{k}$ scaled by the infimum of $\dv{k}$ over $\dvOm$. In general, $\phi$ will depend on the shape of the spatial domain and the thermophysical property distribution.

\item We introduce a  second-order lumped approximation for the domain average QoI $u_\avg$, $\uavgLP$: $\uavgLP$ will  depend on parameters $B$ and also $\phi$. The second-order approximation $\uavgLP$, valid for  \textit{general geometry} and \textit{general thermophysical property distribution}, converges quadraticly in $B$ --- considerably faster convergence in $B$ than the first-order approximation, $\uavgL{1}$.  We also introduce a second-order ``lumped'' approximation for the domain-boundary average QoI $u_\Delta$, $\uDeltaL$; note that, to first order, the temperature variation within $\Omega$ is zero, and hence second-order treatment is required if we wish to obtain a nontrivial approximation for $u_\Delta$.

\noindent We note that our second-order lumped approximation $\uavgLP$ is based on (effectively) Pad\'e approximation, hence slightly different from the Taylor approximation proposed in \cite{ODD} for the homogeneous sphere; see also \cite{Ostro2008}. Our second-order approximation $\uDeltaL$ is also effectively a Pad\'e approximation. In both these approximations, the `P' in the superscript refers to Pad\'e.

\item We present an asymptotic error analysis for all three approximations, $\uavgL{1}$, $\uavgLP$, and $\uDeltaL$,  valid for {general geometry} and {general thermophysical property distribution}. The error analysis relies on $\phi$ as well as two other (similar) quadratic functionals of $\xi$. The asymptotic error estimator for the first-order approximation takes the form $|\uavgL{1}(t) - u_\avg(t)| \leq \phi \,\Bdunk/e + \pcalO(B^2)$; our factor $\phi$ extends the result of \cite{ODD} to general geometry and general thermophysical property distribution. We can also, for the first-order approximation, develop two additional results: we demonstrate that $\uavgL{1}$ is in fact a lower bound for $u_\avg$ for all admissible $t$ and $B$; we provide a \textit{non-asymptotic} error bound for $\uavgL{1}$ which is (less accurate than our asymptotic error estimate but) valid for all $B$. We present companion finite element results to demonstrate the convergence rates and effectivities of the asymptotic error estimators.

\item We provide an explicit engineering interpretation of $\phi$ in terms of a new conduction length scale associated with a \textit{body thermal resistance}. The latter generalizes in a rigorous way the usual definition of thermal resistance --- between two surfaces --- to the notion of thermal resistance between body effective center and body boundary.  Within this interpretation framework, we can easily motivate the particular forms of $\uavgL{1}$, $\uavgLP$, and $\uDeltaL$.

\item We develop closed-form solutions for $\phi$ for the heat transfer canonical spatial domains --- slab, cylinder, and sphere --- with homogeneous thermophysical properties.  (We of course recover $\phi = 3/5$ for the homogeneous sphere, in agreement with \cite{ODD}.) We also develop closed-form exact or asymptotic results for $\phi$ for selected triangles.

 \item We provide a framework for understanding and predicting the dependence of $\phi$ on body composition and spatial domain --- in fact, only body \textit{shape} matters. We provide a useful bound for $\phi$ which requires no information about the spatial distribution of the thermal conductivity, $\dv{k}$. We define a distance between spatial domains which predicts well the stability of $\phi$ with respect to geometric perturbation. We further identify a geometry classification feature: a sufficient condition for spatial domains for which $\phi$ will be large --- and the textbook lumped error criterion potentially misleading. Supporting computational evidence is provided.
\end{enumerate} 

We consider in this paper the case in which the heat transfer coefficient, $\dv{h}$, is independent of time and uniform in space. Furthermore, for the case of natural convection and radiation, the heat transfer coefficient is implicitly linearized about the initial temperature of the body. The more realistic case of temporally-dependent and spatially-dependent heat transfer coefficient is best addressed in the context of a more general discussion which takes as point of departure the ``truth” conjugate heat transfer formulation. In \cite{KandLBandP} [Appendix E] we provide a first theoretical argument, and supporting evidence based on conjugate heat transfer computations, which demonstrates that the results developed in the current paper --- for heat transfer coefficient independent of time and uniform in space --- can in fact  be extended to accommodate small temporal and spatial variations in the local heat transfer coefficient (relative to the average heat transfer coefficient).

We now provide a roadmap for the paper. In Section \ref{sec:form}, we present the dunking initial-value-problem formulation in dimensional and non-dimensional forms; we also introduce the relevant QoIs. In Section \ref{sec:fcnals},  we present and discuss the quadratic functionals of $\xi$ required for our approximations and error estimates. In Section \ref{sec:result}, we introduce our small-Biot lumped approximations for the QoI and associated error estimators;
in Section \ref{sec:numer}, we provide companion finite element results for purposes of assessment. 
In Section \ref{sec:lump}, we motivate, through the resistance formulation, the lumped approximations $\uavgL{1}(t)$, $\uavgLP(t)$, and $\uDeltaL$; we emphasize the physical interpretation of $\phi$ as the correction factor to the usual conduction length scale $\calL$.
Finally, in Section \ref{sec:xi}, we provide several stability and classification results which characterize the dependence of $\phi$ on geometry, and which can serve to identify domains for which $\phi$ --- and hence the error in the first-order lumped approximation --- may be large.

%% file: form.tex
\section{Formulation}\label{sec:form}

The dunking problem can be simply described as follows: a passive solid body at uniform temperature $\dvTi$ is abruptly placed --- dunked --- at time $\dvt = 0$ in an environment, fluid and enclosure, at initial and farfield temperature $\dvTinf$. We recall that a variable with underline is dimensional, whereas a variable without underline is nondimensional. Most of the analysis presented in this work might be extended to more general contexts, for example uniform (Volumetric) heat generation in the body, $\dv{q}^{\textsc{v}}(\dvt)$, or time-dependent farfield temperature, $\dvTinf(\dvt)$. Note that, in general, units are SI, and temperature is in degrees Celsius (or Kelvin).

We first introduce the spatial domain of the body, $\dvOm \subset \RR^d$, and associated boundary $\dvpOm$; we consider spatial dimension $d\in\{1,2,3\}$.  We shall require $\dvOm$ Lipschitz. A point in $\dvOm$ shall be denoted $\dv{x} \equiv (\dv{x}_1,\ldots,\dv{x}_d)$. We choose for our length scale some characteristic dimension of $\dvOm$, $\dv{\ell}$, for example the diameter or the InRadius. We also introduce an intrinsic length scale $\dv{\mathcal{L}}\equiv |\dvOm|/|\dvpOm|$, and an intrinsic associated inverse length scale, 
\begin{align}
\dv{\gamma} \equiv |\dvpOm|/|\dvOm| \label{eq:gammadef_dimensional}\; ; 
\end{align}
here $|\,\cdot |$ denotes measure, and hence, in dimension $d = 3$, $\dv{\gamma}$ is the dimensional surface area ($|\dv{\partial\Omega}|$) to volume ($|\dv{\Omega}|$) ratio. We shall denote time by $\dvt$; we restrict attention to the temporal interval of interest $\dvt \in [0,\dvtf]$.

We now introduce the governing equation for the temperature field within the solid body: the solution $\dvT(\cdot,\dvt)$ to the dunking problem satisfies, for $\dvt\in(0,\dvtf]\,$,
\begin{align}
\dv{\rho c}\, \partial_\dvt \dvT = \dv{\bnabla} \cdot(\dv{k}\, \dv{\bnabla}\, \dvT) \text{ in }\dvOm\,,\label{eq:dim_ivp1}
\end{align}
subject to boundary condition
\begin{align}
\dv{k}\, \partial_\dv{n} \dvT + \dv{h} (\dvT-\dvTinf) = 0 \text{ on }\dvpOm\,,\label{eq:dim_ivp2}
\end{align}
and initial condition
\begin{align}
\dvT(\cdot,0) = \dvTi \text{ in } \dvOm\,. \label{eq:dim_ivp3}
\end{align}
Here, $\partial_\dvt$ and $\partial_\dv{n}$ refer to differentiation with respect to time and domain boundary outward normal direction, respectively.

We now define the nondimensional spatial coordinate as $x \equiv \dv{x}/\dv{\ell}$ and the nondimensional temporal coordinate as $t \equiv \dv{t}/\dv{t}_\text{diff}$, where
\begin{align}
\dv{t}_\text{diff} \equiv \dfrac{\dv{\ell}^2 \,\dashint_{\dv{\Omega}}\dv{\rho c}}{\dv{k}_{\inf}} \, .\label{eq:def_tdiff}
\end{align}
In general, $\dashint$ refers to the average of the integrand over the domain of integration, hence $\dashint_{\underline{\Omega}} \dv{ \rho c} \equiv \dfrac{1}{|\underline{\Omega}|} \int_\Omega \dv{\rho c}$. We further introduce the nondimensional temperature
\begin{align}
u \equiv \dfrac{\dv{T} - \dv{T_\infty}}{\dv{T_{\text{i}}} - \dv{T_\infty}} \, . \label{eq:def_und}
\end{align}
Our nondimensionalization reduces to the standard textbook form in the case in which the thermophysical properties are uniform; notably, $t$ is the Fourier number, $\text{Fo}_{\dv{\ell}}$. We shall implicitly assume that all nondimensional functions take as arguments nondimensional space and/or time.

We next introduce two nondimensional quantities related to thermophysical properties:
\begin{align}
\sigma = \dfrac{\dv{\rho c}}{\dashint_\dvOm \dv{\rho c}} \, ,
\end{align}
\begin{align}
\kappa = \dfrac{\dv{k}}{\dv{k}_{\inf}} \, ,
\end{align}
where $\dv{\rho c}$ and $\dv{k}$ are, respectively, the volumetric specific heat and thermal conductivity, both strictly positive and scalar-valued (generalization to anisotropic heat transfer with symmetric positive definite second-order $\dv{k}$ tensor is not a subject of this work). Note $\dv{\rho}$ is the density and $\dv{c}$ is the mass specific heat; the product $\dv{\rho c}$ is the volumetric specific heat (which is bounded away from zero). Here, $\dv{k}_{\inf}$ (more properly $\dv{k}_{\essinf}$) is the infimum of $\dv{k}$ over $\dvOm$. It follows from our definitions that $\sigma > 0$ and $\dashint_\Omega \sigma = 1$, (and $\kappa > 0$) and $\inf_{\Omega}\kappa = 1$.  Note also that $\sigma = 1$ corresponds to spatially uniform volumetric specific heat, and $\kappa = 1$ corresponds to spatially uniform thermal conductivity.

We may then define the extrinsic Biot number as
\begin{align}
B \equiv \frac{\dv{h}\, \dv{\ell}}{\dv{k}_{\inf}}\,, \label{eq:Bdef}
\end{align}
and the intrinsic Biot number as
\begin{align}
\Bdunk = \dfrac{\dv{h}\, \dv{\calL}}{\dv{k}_{\inf}}\, ; \label{eq:Bidef}
\end{align}
hence $\Bdunk = B\cdot\dv{\calL}/\dv{\ell} = B/\gamma$; here
\begin{align}
    \gamma\equiv \dv{\ell}\,\dv{\gamma}\label{eq:gammadef} \; .
\end{align}
 In the present work, we consider the case in which the heat transfer coefficient, $\dv{h}$, is taken as constant in time and uniform in space. We further assume that the heat transfer coefficient is non-negative, and hence also $B$ is non-negative.

The nondimensional temperature, $u(x,t;B;\sigma,\kappa)$, then satisfies, for any $B \in \Rzp\equiv\{x\in\RR\, |\, x \ge 0\}$,
\begin{align}
\sigma \partial_t u = \bnabla \cdot(\kappa \bnabla u) \text{ in }\Omega\,, \; 0 < t \le \tf\,, \label{eq:ivp1}
\end{align}
subject to boundary condition
\begin{align}
\kappa \partial_n u + B u = 0 \text{ on }\partial\Omega\,, \; 0 < t \le \tf \,, \label{eq:ivp2}
\end{align}
and initial condition
\begin{align}
u = 1 \text{ in } \Omega \text{ at } t = 0 \; ; \label{eq:ivp3}
\end{align}
we shall refer to Equations \eqref{eq:ivp1} -- \eqref{eq:ivp3} as Heat Equation (associated with the dunking problem). Here $\partial_t$ and $\partial_n$ refer to the partial derivatives with respect to time and the outward normal direction, respectively.   We note that, strictly speaking, the strong formulation \eqref{eq:ivp1} -- \eqref{eq:ivp3} requires $\kappa$ sufficiently regular; this difficulty is readily addressed for $\kappa$ piecewise-constant, by a multi-domain strong formulation or, more easily, and more generally, by the weak formulation of the equations \cite{KandLBandP} [Appendix A.2].
We shall in general suppress the parameters $\sigma$ and $\kappa$ unless we are specifically focused on the dependence on these variables. We shall write $u(\cdot,t;B)$ or simply $u(t;B)$ to indicate the spatial nondimensional temperature function for any time $t$ and Biot number $B$ (and, implicitly, prescribed $\sigma$ and $\kappa$); we shall write $u(x,t;B)$ to refer to the evaluation of the spatial nondimensional temperature function at spatial coordinate $x \in \Omega$.

\subsection{Quantities of Interest}\label{sec:QoI}
A Quantity of Interest (QoI) shall refer to a functional of the temperature $u(\cdot,t;B)$. We shall consider in the present work three QoIs, all linear functional outputs: the domain average, denoted $u_\avg$; the boundary average, denoted $u_\paavg$; and the normalized domain average-boundary average difference (more succinctly, domain-boundary average), $u_\Delta$. In our analysis and results, we focus on $u_\avg$ and $u_\Delta$; $u_\paavg$ is a necessary intermediate definition. The domain average QoI is defined as
\begin{align}
u_\avg(t;B) \equiv  \dashint_{\Omega} \sigma u(\cdot,t;B) \; ;\label{eq:uavg}
\end{align}
the boundary average QoI is defined as
\begin{align}
u_\paavg(t;B) \equiv  \dashint_{\partial\Omega} u(\cdot,t;B) \; ;\label{eq:upaavg}
\end{align}
and the (normalized) domain-boundary average QoI is defined as
\begin{align}
u_\Delta(t;B) \equiv \dfrac{u_\avg(t;B) - u_\paavg(t;B)}{u_\avg(t;B)} \; .\label{eq:udelta}
\end{align}
These QoIs are defined for almost all times $t$ in our interval $[0,\tf]$. Note it follows from the maximum principle \cite{Evans} that $u_\avg$ and $u_\paavg$ take on values in the range $[0,1]$. We choose to include $\sigma$ in the integrand of \eqref{eq:uavg} --- hence a weighted average --- to reflect the physical quantity (and associated conservation properties) represented by this QoI: $u_\avg$ is the nondimensional volume-specific thermal energy.

We note that the QoI, in particular the domain average QoI, can serve in either forward or inverse mode. In forward mode, given $B$, we wish to evaluate $u_\avg(t;B), 0 \le t \le \tf$. In inverse mode, given a target domain average QoI value, $u^*_\avg$, we wish to evaluate $\tau_\avg(u_\avg^*;B)$ such that
\begin{align}
u_\avg(\tau_\avg(u_\avg^*;B);B) = u_\avg^* \; .
\end{align}
For example, if $u_\avg^* = \exp(-1)$, then we might interpret $\tau_\avg$ as the time constant of the system. In some cases, $\tau_\avg$ can be a more meaningful and sensitive output than $u_\avg$. The error estimators for $\tau_\avg$ follow directly from the error estimators for $u_\avg$; we shall not further consider the inverse mode in this paper.

%% file: fcnals.tex
\section{Sensitivity Formulation}\label{sec:fcnals}

We first  present the sensitivity equation for our field $\xi$. We then proceed to the quadratic functionals of $\xi$ required for the subsequent development.

\subsection{Sensitivity Equation}

We  introduce here the equation for the sensitivity derivative $\xi$, a key field in the synthesis and analysis of the small-Biot approximations. The field $\xi$ satisfies
\begin{align}
- \bnabla \cdot (\kappa \bnabla \xi) = \OmMeas^{-1/2}\gamma \sigma \text{ in } \Omega \; \label{eq:xi1}
\end{align}
subject to boundary condition
\begin{align}
\kappa \partial_n \xi = -\OmMeas^{-1/2} \text{ on } \partial\Omega \;  \label{eq:xi2}
\end{align}
and zero-mean condition
\begin{align}
\int_\Omega \sigma \xi = 0 \; . \label{eq:xi3}
\end{align}
We recall that $|\Omega| = \int_\Omega 1$ and (for future reference) $|\partial\Omega| = \int_{\partial\Omega} 1$; we also recall that $\partial_n$ denotes the derivative with respect to the outward normal on $\partial\Omega$. 

We note that \eqref{eq:xi1} -- \eqref{eq:xi3}, a constrained elliptic PDE, is solvable --- since $\sigma$ is of average unity over $\Omega$ --- and furthermore admits a unique solution --- thanks to the condition \eqref{eq:xi3}. We refer the reader to \cite{KandLBandP} [Appendix C.7] for additional details on the derivation and analysis of \eqref{eq:xi1} -- \eqref{eq:xi3} (note that  $\xi$ is denoted by $\psipz$ in \cite{KandLBandP}).

\subsection{Quadratic Functional Outputs}

\subsubsection{Definition}

We introduce here three quadratic functional outputs --- $\phi$, $\chi$, and $\Upsilon$ --- of $\xi$; these shall play an important role in the analysis. We have already discussed the ubiquitous role of $\phi$; the other two functional outputs, $\chi$ and $\Upsilon$, appear in the error estimators for our second-order (Pad\'e) approximations.

We provide the definitions:
\begin{align}
\phi&\equiv\int_\Omega\kappa\, \nabla\xi\cdot\nabla\xi\; ,\label{eq:phidef}\\
\chi &\equiv \int_{\partial\Omega}\xi^2\; , \label{eq:chidef}\\
\Upsilon&\equiv\int_\Omega \sigma\,\xi^2\; ;\label{eq:Upsilondef}
\end{align}
for future reference we also note, from \eqref{eq:xi1} -- \eqref{eq:xi3}, that
\begin{align}
\int_\Omega\kappa\, \nabla\xi\cdot\nabla\xi = -|\Omega|^{-1/2}|\partial\Omega| \dashint_{\partial\Omega}\xi \; , \label{eq:kappa_alt}
\end{align}
which thus serves as an alternative (equivalent) expression for $\phi$. A physical interpretation of $\phi$ --- the most important of the trio $\phi, \Upsilon$, and $\chi$ --- is discussed in Section \ref{sec:lump}; however, we define already here
\begin{align}
\Bdunkp\equiv \frac{\phi B}{\gamma},
\end{align}
a $\phi$-corrected Biot number which shall appear frequently in the  subsequent analysis. We emphasize that all three quadratic functional outputs may be inexpensively evaluated in terms of $\xi$,  and thus the principal computational task is the solution of \eqref{eq:xi1} -- \eqref{eq:xi3} for $\xi$; we discuss computational issues and cost, in context, as we proceed.

This PDE formulation to obtain $\phi$ (and $\chi$ and $\Upsilon$) is the foundation of our framework. Other formulations for quadratic functional outputs are also possible, in particular direct Taylor-series eigenvalue expansion; an example of the latter is developed in \cite{ODD} for $\phi$ for the special case of a spherical body with homogeneous properties ($\sigma = \kappa = 1$). We claim that the PDE formulation \eqref{eq:xi1} -- \eqref{eq:xi3}, \eqref{eq:phidef} -- \eqref{eq:Upsilondef} offers critical advantages: more systematic theoretical analysis and also numerical computation; more ready accommodation of non-uniform thermophysical properties ($\sigma \neq 1, \kappa \neq 1$); more flexible treatment of geometry and in particular general spatial domains $\Omega$. Indeed, even closed-form solutions (in simple geometries) are much more readily revealed by the PDE formulation  \eqref{eq:xi1} -- \eqref{eq:xi3}, \eqref{eq:phidef} -- \eqref{eq:Upsilondef} --- a simple boundary value problem rather than a nonlinear characteristic equation; we present several useful closed-form solutions in the next section.

\subsubsection{Properties}

We summarize here several important properties for our quadratic functional outputs $\phi(\Omega;\sigma,\kappa)$, $\chi(\Omega;\sigma,\kappa)$, and $\Upsilon(\Omega;\sigma,\kappa)$. We do place greater emphasis on $\phi$ since only $\phi$ is implicated in the first-order (approximation and) error estimator. We enumerate the properties for later reference:

\begin{itemize}
\item[P1] Spatial Scale Invariance: $\phi(\Omega_1;\ldots) = \phi(\Omega_2;\ldots)$, for any domain $\Omega_2$ which is a translation, rotation, and dilation of $\Omega_1$. In contexts in which we wish to emphasize spatial scale invariance, we shall in arguments for $\phi$ refer to $\Omega_{\text{ref}}$ (rather than $\Omega$), where $\Omega_{\text{ref}}$ is reference shape for a family of similar domains; recall that $\Omega_{\text{ref}} \subset \RR^d$ for $d \in \{1,2,3\}$.
\item[P2] Positivity: $\phi(\Omega_{\text{ref}}\,; \sigma,\kappa) > 0$, $\chi(\Omega_{\text{ref}}\,; \sigma,\kappa) > 0$, $\Upsilon(\Omega_{\text{ref}}\,; \sigma,\kappa) > 0$ for any given admissible $\sigma$ and $\kappa$.
\item[P3] Monotonicity\footnote{Also in \cite{KandLBandP} [Proposition 5.3] we develop some bounds for the $\sigma$-dependence of $\phi$, however the results are not sufficiently general to warrant inclusion here.} in $\kappa$:
$\phi(\Omega; \sigma,\kappa_1) \ge \phi(\Omega;\sigma,\kappa_2)$,
for $\kappa_1$ and $\kappa_2$ such that $\kappa_1 \le \kappa_2, \forall x \in \Omega$. Corollary: $\phi(\Omega;\sigma,1)$ is an upper bound for $\phi(\Omega;\sigma,\kappa)$ for any admissible $\kappa$.
\item[P4]  Construction of $\phi$ for Tensorized Domains: Given $m \subset \{2,3\}$,  a domain $\Omega^{[m-1]} \subset \RR^{m-1}$, and $L_{m} > 0$, define the extruded (right-cylinder) domain $\Omega \subset \RR^m$ as $\Omega \equiv \Omega^{[m-1]} \times (0,L_{m})$.  Then $\phi(\Omega;1,1) = \phi(\Omega^{[m-1]};1,1) + \dfrac{1}{3}$.
\item[P5] Spatial Scale Invariance: $(\gamma_1)\chi(\Omega_1;\ldots) = (\gamma_2)\chi(\Omega_2;\ldots)$, for any domain $\Omega_2$ which is a translation, rotation, and dilation of $\Omega_1$; here $\gamma_1$ and $\gamma_2$ are respectively the values of $\gamma$ associated with $\Omega_1$ and $\Omega_2$. In contexts in which we wish to emphasize spatial scale invariance, we shall in arguments for $\chi$ refer to $\Omega_{\text{ref}}$ (rather than $\Omega$), where $\Omega_{\text{ref}}$ is reference shape for a family of similar domains; recall that $\Omega_{\text{ref}} \subset \RR^d$ for $d \in \{1,2,3\}$.
\item[P6]  Construction of $\chi$ for Tensorized Domains: Given $m \subset \{2,3\}$,  a domain $\Omega^{[m-1]} \subset \RR^{m-1}$, and $L_{m} > 0$, define the extruded (right-cylinder) domain $\Omega \subset \RR^m$ as $\Omega \equiv \Omega^{[m-1]} \times (0,L_{m})$.  Then $\chi(\Omega;1,1) = \chi(\Omega^{[m-1]};1,1) + \chi(\Omega^{[1]};1,1) + \frac{2}{L_m}\Upsilon(\Omega^{[m-1]};1,1)+\gamma_{[m-1]}\Upsilon(\Omega^{[1]};1,1)$; $\gamma_{[m-1]}$ is the value of $\gamma$ associated with domain $\Omega^{[m-1]}$.
\item[P7] Spatial Scale Invariance: $(\gamma_1)^2\Upsilon(\Omega_1;\ldots) = (\gamma_2)^2\Upsilon(\Omega_2;\ldots)$, for any domain $\Omega_2$ which is a translation, rotation, and dilation of $\Omega_1$; here $\gamma_1$ and $\gamma_2$ are respectively the values of $\gamma$ associated with $\Omega_1$ and $\Omega_2$. In contexts in which we wish to emphasize spatial scale invariance, we shall in arguments for $\Upsilon$ refer to $\Omega_{\text{ref}}$ (rather than $\Omega$), where $\Omega_{\text{ref}}$ is reference shape for a family of similar domains; recall that $\Omega_{\text{ref}} \subset \RR^d$ for $d \in \{1,2,3\}$.
\item[P8]  Construction of $\Upsilon$ for Tensorized Domains: Given $m \subset \{2,3\}$,  a domain $\Omega^{[m-1]} \subset \RR^{m-1}$, and $L_{m} > 0$, define the extruded (right-cylinder) domain $\Omega \subset \RR^m$ as $\Omega \equiv \Omega^{[m-1]} \times (0,L_{m})$.  Then $\Upsilon(\Omega;1,1) = \Upsilon(\Omega^{[m-1]};1,1) + \Upsilon(\Omega^{[1]};1,1)$.
\end{itemize}
These properties are all proven in \cite{KandLBandP} [in order, Proposition 5.2, Proposition 5.1, Proposition 5.3, Proposition 5.7, Proposition 6.8, Proposition 6.9, Proposition 6.8, and Proposition 6.9]. We shall see later how these properties, in particular P3 and P4, can serve to develop closed-form solutions or bounds for $\phi$.

\subsubsection{Closed-Form Solutions}\label{sec:cfs}

We present in Table \ref{table:canonicalchi} values for our three quadratic functional outputs for uniform thermophysical properties ($\sigma = \kappa = 1$) for several canonical domains. We recover, in agreement with the earlier work \cite{ODD},  the value $\phi  = 3/5$ for the homogeneous sphere; the details of our derivation --- and the very simple polynomial form for $\xi$ --- are provided in \cite{KandLBandP} [Table 7]. The real utility of the values for these canonical domains is in conjunction with the properties P3, P4, P6, and P8 and the stability conjecture presented in Section \ref{sec:xi}.

\input{Tables/canonfeval}

 We first illustrate the application of P3. Let us consider, for simplicity, a sphere with $\sigma = 1$ uniform and $\kappa \ge 1$ non-uniform. We recall from the introduction that our asymptotic estimate for the error $|u_\avg - \uavgL{1}|$ is $\phi B/(\gamma e)$. In general, $\phi$ will depend on $\kappa$; however, thanks to P3, we may bound $\phi(\Omega;1,\kappa)$ by $\phi(\Omega;1,1) = 3/5$ from Table \ref{table:canonicalchi}.
 
 We next illustrate tensorization. In all cases, we consider exclusively $\sigma = \kappa = 1$, and we focus on $\phi$.  As a first instance, we tensorize the interval with the interval to obtain $\phi = 2/3$ for a rectangle; we can subsequently tensorize our rectangle with the interval to obtain $\phi = 1$ for a rectangular parallelepiped (slab). Note it follows from scale invariance of the building blocks that $\phi$ does not depend on the dimensions (aspect ratio) of the rectangle or parallelepiped. As a second instance, we tensorize the disk with the interval to obtain $\phi = 1/2 + 1/3 = 5/6$ for a right circular cylinder; again, $\phi$ does not depend on radius or axial length.  As a third instance, we tensorize $\Omega_{\text{right triangle}}(W)$ with the interval to obtain a triangular prism: $\phi$ will depend only on $W$.
 In actual practice, only three-dimensional spatial domains are of interest; however, thanks to the tensorization property, we can (for extruded domains) restrict attention to two-dimensional spatial domains --- easier computations and visualizations.

%% file: Tables/canonfeval.tex
\begin{table}[H]
\begin{center}
\caption{Values of the $\gamma$-scaled (scale-invariant) quadratic functional outputs for uniform thermophysical properties ($\kappa = 1, \sigma = 1$) and several canonical geometries. The domain $\Omega_\text{right triangle}(W)$ is defined by vertices $(0,0)$, $(W,0)$, and $(0,1)$.}  \label{table:canonicalchi}
\begin{tabular}{l|c|c|c}
\multicolumn{1}{c|}{$\Omega$} & $\phi$ & $\gamma\,\chi$ & $\gamma^2\,\Upsilon$\\
\hline\hline
$\Omega_\text{interval}$  & 1/3 & 1/9 & 1/45\\
$\Omega_\text{disk}$  & 1/2 & 1/4 & 1/12\\
$\Omega_\text{sphere}$ & 3/5 & 9/25 & 27/175\\
$\Omega_\text{isosceles right triangle}$ & 4/3 & $\frac{4}{5}(3+2\sqrt{2})$ & $\frac{4}{15}(3+2\sqrt{2})$\\
$\Omega_\text{equilateral triangle}$ & 1 & 9/5 & 3/5\\
$\Omega_\text{right triangle}(W)$ & $\sim \frac{2}{3}W^{-2} $ & $\sim \frac{28}{15} W^{-4}$ & $ \sim \frac{28}{45} W^{-4}$ 
\end{tabular}
\end{center}
\end{table}

%% file: result.tex
\section{Approximation and Error Estimation}\label{sec:result}

Here we introduce the small-Biot approximations and associated asymptotic error estimators for the QoIs \eqref{eq:uavg} -- \eqref{eq:udelta}. We focus here on the mathematical descriptions of these quantities; we consider the physical interpretations in Section \ref{sec:lump}.

\subsection{Small-Biot Approximations}

We first  introduce the small-Biot approximations to the QoIs \eqref{eq:uavg} -- \eqref{eq:udelta}, introduced in Section \ref{sec:QoI}. Associated with $u_\avg(t;B)$, we consider two approximations: the first is the first-order ``classical'' lumped approximation, $\uavgL{1}(t;B)$, prevalent in engineering applications; the second is a second-order lumped approximation, $\uavgLP(t;B)$. Finally, we introduce the second-order approximation to $u_\Delta(t;B)$, $\uDeltaL(t;B)$; we note that the first-order approximation to $u_\Delta(t;B)$ is simply zero.

We now define our lumped approximations:
\begin{align}
\uavgL{1}(t;B) &\equiv \exp(-B\gamma t) \; , \label{eq:u1avg}\\
\uavgLP(t;B) &\equiv \exp\left(-\dfrac{B\gamma t}{1 + \phi B/\gamma}\right) \; , \label{eq:u2avg}\\
\uDeltaL(B) &\equiv \dfrac{\phi B}{\gamma}\left(1 + \frac{\phi B}{\gamma}\right)^{-1} \; . \label{eq:u2Delta}
\end{align}
For future reference, we also define from \eqref{eq:u1avg}  the first-order approximation for the ``equilibration time constant'' for the domain average QoI:
\begin{align}
\taueq \equiv  \dfrac{1}{B\gamma} \; . \label{eq:taueqdef}
\end{align}
Note that $\uDeltaL(B)$ is independent of time. We can deduce from Equations \eqref{eq:u1avg} -- \eqref{eq:u2Delta} and the definition of $\uDeltaL$, \eqref{eq:udelta}, an approximation for the boundary average temperature, $u_\paavg$: $\tilde{u}^{2\text{P}}_{\paavg} (t;B) \equiv \uavgL{1}(t;B) (1+\uDeltaL(B))$. From these definitions, we observe that the functional output, $\phi$ \eqref{eq:phidef}, plays a key role in the second-order approximations, $\uavgLP(t;B)$ and $\uDeltaL(B)$. In the following section, we will establish that $\phi$ is also the key new ingredient in the error analysis of $\uavgL{1}(t;B)$.

We note that Gockenbach and Schmidtke propose in \cite{ODD} a second-order {\em Taylor}-series approximation for the particular case of a homogeneous sphere. In contrast, our second-order approximation (valid for general thermophysical properties and spatial domains) is in fact a \textit{Pad\'e} approximation  --- hence the P in the superscript $\uavgLP$. The Taylor and Pad\'e approximations will yield the same (quadratic) asymptotic convergence rate with $B$ as $B \to 0$; however, as is often the case, the Pad\'e approximation $\uavgLP$ yields meaningful results over a larger range of $B$. In a similar fashion, the denominator in \eqref{eq:u2Delta} ensures that (the Pad\'e approximation) $\uDeltaL$ is meaningful --- and always less than unity --- over a larger range of $B$. In Section \ref{sec:approx_phys} we shall demonstrate that our Pad\'e approximations naturally arise from resistance considerations.

\subsection{Error Analysis for $\uavgL{1}$}\label{sec:ea}

In this section, we establish the error estimate for the $\uavgL{1}$ approximation \eqref{eq:u1avg}. We remind the reader that the well-established criterion for lumped approximation is $\Bdunk$ small: true, but not quantitative. Now, equipped with the general PDE analysis framework, we can estimate the approximation error in $\uavgL{1}$ for any given spatial domain and thermophysical property distribution. 

We shall define the error in $\uavgL{1}$ as
\begin{align}
    \eavgL{1}(t;B) \equiv \max_{s\in[0,t]} |u_\avg(s;B) - \uavgL{1}(s;B)| \; . \label{eq:eavgLdef}
\end{align}
We establish two types of bounds for $\eavgL{1}$: asymptotic and non-asymptotic. The asymptotic bound is applicable only in the small-Biot regime, in particular as $B \rightarrow 0$; the non-asymptotic bounds are valid for all $B > 0$. 

We first establish the asymptotic bound for $e^1_\avg(t;B)$ \cite{KandLBandP} [Proposition 6.5]:
\begin{align}
e^1_\avg(t;B) \le \frac{\phi B}{\gamma}\exp(-1) + \pcalO(B^2) \text{ as } B\rightarrow 0. \label{eq:uavg1abound}
\end{align}
We note that the leading-order error term
\begin{align}
\eavgLa{1}\equiv\frac{\phi B}{\gamma}\exp(-1)\label{eq:uavg1abounddef}
\end{align}
is proportional to $\phi B/\gamma$, the $\phi$-corrected intrinsic Biot number $\Bdunkp$. 

Next we introduce the non-asymptotic error bound for $\uavgL{1}$ \cite{KandLBandP} [Proposition 6.6]:
\begin{align}
\eavgL{1}(t;B) \,\le\, \frac{1}{2}\sqrt{\frac{\phi B}{\gamma}}\,,\ \forall\, t \in [0,\tf]\;. \label{eq:uavg1bound}
\end{align}
We note that for this non-asymptotic result, the error bound
\begin{align}
\eavgLUB{1} \equiv \frac{1}{2} \biggl(\underbrace{\frac{\phi B}{\gamma}}_{\Bdunkp}\biggr)^{1/2}\;
\end{align}
still depends on the $\phi$-corrected intrinsic Biot number $\Bdunkp$. The non-asymptotic bound \eqref{eq:uavg1bound} eliminates any uncertainty in the error estimate. However, it is clear that the non-asymptotic bound is not as sharp as the asymptotic bound \eqref{eq:uavg1abound} as $B \rightarrow 0$.  We refer the reader to \cite{KandLBandP} for derivation of these error bounds [Appendix C.17 and Appendix C.18].

We can further prove \cite{KandLBandP} [Proposition 6.3]
\begin{align}
\uavgL{1}(t;B) \le u_\avg(t;B), \forall t \in [0,\tf] \;; \label{eq:LBallB}
\end{align}
notably, the relation \eqref{eq:LBallB} does not require evaluation of $\phi$. The physical interpretation of this bound is standard: the lumped assumption of temperature uniformity implicitly neglects the body (internal) conduction resistance, hence overestimates the heat transfer rate, hence underestimates the equilibration time; we elaborate on this argument in Section \ref{sec:lump} and Section \ref{sec:xi}. 
It is clear that, for $B \gg 1$, $\uavgL{1}$ will severely underestimate $u_\avg$. However, for $B$ order unity (and certainly for small $B$), the approximation $\uavgL{1}$ can prove useful --- in particular in engineering contexts in which a lower bound for $u_\avg$ can provide useful guidance to a design or decision process.

Finally, several comments on utility. We first discuss the relevance to assessment and error control. We observe that, if $\phi \gg 1$, then --- if the error estimator is sharp --- the error in the first-order lumped approximation will be much larger than predicted by the usual textbook criterion; in particular, even if $\Bdunk$ is small, $\Bdunkp \equiv \phi \Bdunk$ may not be small. We further note from Table \ref{table:canonicalchi} that there indeed exist spatial domains for which $\phi \gg 1$. We shall conclude this argument in the next section: a large error estimate, in general, does not imply a large error; however, in Section \ref{sec:numer} we shall confirm (for a particular problem)  that the error estimate is reasonably sharp, hence representative of the actual error.

We next discuss computational cost. We recall that our first-order lumped approximation, $\uavgL{1}$, is a surrogate for $u_\avg$, the domain average QoI associated with the solution $u$ to Heat Equation (hereafter abbreviated HE), \eqref{eq:ivp1}~--~\eqref{eq:ivp3}. The advantages of $\uavgL{1}$ relative to $u_\avg$ are threefold: transparency — explicit parameter dependence; flexibility — ready incorporation into design or optimization or system analysis; rapid calculation — just a few standard elementary function evaluations. However, if we wish to certify $\uavgL{1}$ with our error estimator $e^{\text{1 asymp}}_\avg$, we  must re-introduce a PDE, in particular \eqref{eq:xi1} -- \eqref{eq:xi3}. 

From a computational perspective we would, of course, prefer to avoid the cost associated with solving a PDE. However, there are mitigating considerations. First, in some cases we can avoid numerical solution of the PDE for $\xi$: we can appeal to the closed-form solutions (and related bounds) of Section \ref{sec:fcnals}; we will further expand the reach of the closed-form dictionary  --- to admit geometry perturbations --- based on the stability arguments to be introduced in Section \ref{sec:xi}. Second, our error estimator requires numerical solution of an {\em elliptic} PDE, whereas $u_\avg$ requires numerical solution of a {\em parabolic} PDE, in particular HE; in general, the former is considerably less expensive than the latter.  Third, it is important to note that, although our error estimator may indeed increase the computational cost --- and perhaps software effort\footnote{The software issues may be addressed by microservices; see \cite{mathworks} for a discussion of microservices in the heat transfer context.} --- associated with our (certified) lumped approximation, we retain key advantages of the lumped approximation $\uavgL{1}$, in particular related to transparency and interpretation. 

\subsection{Error Analysis for $\uavgLP$}\label{sec:ea2}

Proceeding in a similar fashion, we now provide error estimates for our second-order lumped Pad\'e approximation, $\uavgLP$ \eqref{eq:u2avg}. We first define
\begin{align}
\eavgLP(t;B) \equiv \max_{s \in [0,t]} |u_\avg(s;B) - \uavgLP(s;B)|.
\end{align}
Then the asymptotic error satisfies \cite{KandLBandP} [Proposition 6.5]
\begin{align}
\eavgLP(t;B) \le \eavgL{2\text{P\,asymp}} + \pcalO(B^3) \text{ as } B \rightarrow 0,
\end{align}
where
\begin{align}
\eavgL{2\text{P\,asymp}}\equiv\left(\frac{|\gamma\chi-\gamma^2\Upsilon-\phi^2|}{e\gamma^2}+\Upsilon\right)B^2,\label{eq:uavg2abound}
\end{align}
which may also be expressed in terms of scale-invariant terms and $\Bdunk$ as
\begin{align}
\eavgL{2\text{P\,asymp}}=\left(|\gamma\chi-\gamma^2\Upsilon-\phi^2|\exp(-1)+\gamma^2\Upsilon\right)\Bdunk^2.\label{eq:sinv2p}
\end{align}
Thus, once $\xi$ is known --- and hence also (at negligible additional cost) $\phi$, $\chi$, and $\Upsilon$ are known --- we can evaluate both the second-order Pad\'e approximation, $\uavgLP$, and associated error estimate, $\eavgL{2\text{P\,asymp}}$.

\subsection{Error Analysis for $\uDeltaL$}\label{sec:ea3}
We now establish an error bound for the domain-boundary average QoI, $\uDeltaL$ of  \eqref{eq:u2Delta}. First we define the error quantity for this approximation:
\begin{align}
\eDeltaLr(t_0,\tf;B) \equiv \max_{s \in [t_0,\tf]} \dfrac{|u_\Delta(s;B) - \uDeltaL(B)|}{\uDeltaL(B)} \; , \label{eq:eDeltaLrdef}
\end{align}
for some $t_0 \in (0,\tf)$.
Next, we introduce scale-invariant constants in terms of our quadratic functionals of $\xi$,
\begin{align}
C_0&\equiv \frac{\gamma^2\Upsilon}{\phi}\exp(-1) \;,\\
C_1&\equiv \frac{|\gamma\chi-\gamma^2\Upsilon - \phi^2|}{\phi}\; .
\end{align}
We may then define the error bound:
\begin{align}
\eDeltaLr(t_0,\tf;B) \le \left(\frac{C_0}{B\gamma t_0} + C_1\right) \underbrace{(B/\gamma)}_{\Bdunk} +\, \pcalO(B) \text{ as }B\rightarrow 0\; ; \label{eq:eDeltarelbound}
\end{align}
the terms reflected in $\pcalO(B)$ are independent of $\tf$ \cite{KandLBandP} [Proposition 6.10]. We shall refer to $(C_0/(B\gamma t_0))\Bdunk$ as the transient contribution to the asymptotic relative error estimate, and $C_1 \Bdunk$ as the long-time contribution to the asymptotic relative error estimate. We emphasize that $\eDeltaLr(t_0,t;B)$ is the maximum error in the lumped approximation {\em relative} to the lumped approximation for $u_\Delta(t;B)$, $\uDeltaL(B)$, over the time interval $[t_0,\tf]$. In the limit $B \rightarrow 0$, although $\uDeltaL(t;B)$ is $\pcalO(B)$, the {\em absolute} error in $\uDeltaL(t;B)$ is order $\pcalO(B^2)$ --- which leads to the first-order convergence result for the relative error, \eqref{eq:eDeltarelbound}. In short, we can accurately discriminate $u_\Delta(t;B)$ as a function of $B$.

We can not hope to capture the short-time evolution of the temperature field on the boundary within our this second-order ``lumped'' framework. The cut-off $t_0$ restricts attention to times away from $t = 0$; as $t_0$ decreases, the error in $\uDeltaL(t;B)$ increases commensurately. However, we note that our error estimate \eqref{eq:eDeltarelbound} is valid for $t \ge t_0$, and hence for $t_0$ ``small'' (say, $t_0 = 0.2 \taueq$), the error estimate is relevant over most of the time interval of interest, $(0,\tf]$ (say for $\tf\ge \taueq)$). We recall that $\taueq \equiv1/(B\gamma)$, defined in \eqref{eq:taueqdef}, is the first-order approximation for the domain average QoI equilibration time constant.

%% file: numer.tex
\section{Numerical Results}\label{sec:numer}

We first introduce two triangular domains: SART-1 and SART-2, where SART is an acronym for Small Angle Right Triangle. Here SART-1 (and later variants) refers to a right triangle with legs in the ratio 4:1 (and then corresponding hypotenuse), hence $\Omega_{\text{right triangle}}(W = 1/4)$; SART-2 refers to a right triangle with legs in the ratio 16:1 (and then corresponding hypotenuse), hence $\Omega_{\text{right triangle}}(W = 1/16)$. In the current section, focused on assessment of the lumped approximations and associated error estimators, we consider the more extreme SART-2 --- in particular to emphasize error amplification by $\phi$; in Section \ref{sec:xi}, focused on the dependence of $\phi$ on domain perturbation, we consider the less extreme, hence more easily visualized, SART-1.

We now proceed to assess our lumped approximations and associated error estimators for domain SART-2 by direct comparison with finite element numerical results for the dunking parabolic PDE, Heat Equation. We consider the case of uniform thermophysical properties, $\sigma = \kappa = 1$.
The numerical results in this section are based on a finite-difference in time finite-element in space discretization with uniform refinement and associated {\it a posteriori} error estimators; we exploit the latter to ensure that the discretization error will always be sufficiently small compared to the actual errors we wish to investigate --- the errors associated with our lumped approximations.

As described in Section \ref{sec:result}, our lumped approximations and error estimators will depend (for our given domain) on $B$, which we shall vary, as well as the our quadratic functionals of $\xi$. We provide the values for the latter: $\phi = 161$, $\gamma\chi = 1.21\cdot 10^5$; $(\gamma)^2\Upsilon  = 4.02\cdot 10^4$; these quantities are computed by the adaptive finite element procedure summarized in Section \ref{sec:xi}. Note that the value $\phi = 161$ is quite close to the asymptotic result, $\phi \sim (2/3)W^{-2} = 170.67$, presented in the last row of Table \ref{table:canonicalchi} of Section \ref{sec:cfs}.

\subsection{Domain Average Approximation $\uavgL{1}$}

We show in Table \ref{tab:NumRes2a} the true error, the  asymptotic error estimate, and the non-asymptotic error bound for the first-order lumped approximation of $u_\text{avg}$, $\uavgL{1}$,  as a function of Biot number $B$ (and $\Bdunk$) for domain SART-2. We choose for our final time $\tf = 2 \taueq$; it is shown in \cite{KandLBandP} [Proposition 6.5] that (for $\tf \ge \taueq$) the error in $\uavgL{1}(t;B)$ is a maximum at time $t = \taueq$ as $B \rightarrow 0$; hence our results are in fact valid for any $\tf > \taueq$ (and conservative for $\tf < \taueq$). 

We make the following observations from Table \ref{tab:NumRes2a}. We confirm that the error $\eavgL{1}$, tends to zero linearly in $B$: for example $\eavgL{1}(\tf)$ for $B=2\cdot 10^{-3}$ divided by $\eavgL{1}(\tf)$ for $B=1\cdot 10^{-3}$ is 1.98, quite close to the theoretical value of 2 as $B\to 0$; the first-order approximation is indeed first-order.
We confirm that our error estimator, $e_\text{avg}^\text{1 asymp}$, is indeed an upper bound as $B \to 0$; we further observe that $e_\text{avg}^\text{1 asymp}$ is asymptotically exact --- $e_\text{avg}^\text{1 asymp}/\eavgL{1} \rightarrow 1$ (from above) as $B \to 0$. We confirm that our non-asymptotic error bound, $e_\text{avg}^\text{1 UB}$, is indeed an upper bound for the true error for all $B$ (presented) --- but clearly not as sharp as $e_\text{avg}^\text{1 asymp}$ as $B \rightarrow 0$. Finally, we note that if $\phi$ is replaced by an upper bound for $\phi$ --- possible in some circumstances, as described in Section \ref{sec:xi} --- then we would retain both our asymptotic and non-asymptotic bound properties, though the error estimators would of course be less sharp.

We emphasize that, for (say) $B = 0.01$, our asymptotic error estimate is very sharp and also an upper bound. Had we simply set $\phi = 1$ in our error estimators --- as might be suggested by the usual rule of thumb that the error in the (first-order) lumped approximation should be on the order of $\Bdunk$ --- we would have {\em under}estimated the error by roughly two orders of magnitude: our error estimator is $\phi \Bdunk$ (for $\phi = 161$), whereas the rule of thumb would suggest simply $\Bdunk$. The latter is imprecise, and even worse, very optimistic; thus, the importance of $\phi$. In Section \ref{sec:lump} and Section \ref{sec:xi} we will provide several physical interpretations for $\phi$. We will then be positioned in Section \ref{sec:xi} to propose strategies for estimation of $\phi$; of particular importance are geometric features which lead to large values of $\phi$. 

Before further investigation of $\phi$, however, we shall complete our assessment of our lumped approximations --- it remains to study $\uavgLP$ and $\uDeltaL$ --- and associated error estimators.

\input{Tables/sart2a.tex}

\subsection{Domain Average Approximation $\uavgLP$}

We  present  in Table \ref{tab:EavgL2_sart2} the numerical results for the error $\eavgLP(\tf)$ for the same problem parameters (including $\tf = 2\taueq$) studied in the previous section: we observe second-order convergence and asymptotic bounds, as expected. Note, and in contrast to the first-order case, the effectivity of the error estimator, here $e_\text{avg}^\text{2P asymp}/\eavgLP$, does not approach 1 as $B \to 0$. However, and more importantly, the second-order approximation does deliver much improved accuracy relative to the first-order approximation for small $B$.  In theory this improvement in accuracy comes at very little cost, since once we have obtained $\xi$ (and $\phi$), $\chi$ and $\Upsilon$ are effectively ``for free." However, in practice, the the merit of the second-order approximation is less clear: for the first-order approximation and error estimator we need only $B$ and $\phi$ {\em or} an upper bound estimate for $\phi$, whereas for the second-order approximation we require --- if we wish to retain second-order convergence --- precise values of $\phi$, $\chi$, and $\Upsilon$; we discuss this point further in Section \ref{sec:xi}. 

\input{Tables/sart2b.tex}

\subsection{Domain-Boundary Average Approximation $\uDeltaL$}

We consider the same problem studied in the previous two sections. We recall that $\tf = 2\taueq$. We must now also choose $t_0$: we take $t_0 = 0.2 \taueq$, and hence our error estimator is defined for most of the time interval of interest. We present in Table \ref{tab:EDeltaL_sart2} the relative error and associated asymptotic relative error estimator for the second-order approximation of $u_\Delta$, $\uDeltaL$, as a function of $B$ (and $\Bdunk$). The third column is the true relative error, $\eDeltaLr(t_0,\tf;B)$. The fourth and fifth columns derive from the asymptotic relative error estimator provided in \eqref{eq:eDeltarelbound}: the fourth column is the long-time contribution to the asymptotic relative error estimator, and the fifth column is the full asymptotic relative error estimator (expressed as the sum of the transient and long-time contributions).
We observe, as predicted from the theory, first-order convergence in $B$ and asymptotic bounds. We note that the asymptotic relative error estimator is not overly sharp as $B \to 0$, which, from the fourth and fifth columns of Table \ref{tab:EDeltaL_sart2}, we can attribute to the transient contribution; the long-time contribution alone provides a much sharper estimator, but in theory (and in practice) may not constitute an upper bound.

We emphasize that the classical lumped approximation --- a first-order approximation --- provides no information for $u_\Delta$; more precisely, to first order, $u_\Delta = 0$ for all time.  In contrast, the second-order ``lumped'' approximation\footnote{We recall that, although the relative error in $\uDeltaL$ is $\pcalO(B)$, the absolute error in $\uDeltaL$ is $\pcalO(B^2)$ --- hence the label ``second-order."}  provides valuable information about the temperature variation within the body, and hence also the temperature on the boundary. We again observe the central role of $\phi$, in this instance appearing not only in the error estimator (through $C_0$ and $C_1$) but also directly in our approximation: $\uDeltaL(B) \equiv \frac{\phi B}{\gamma}\left(1 + \frac{\phi B}{\gamma}\right)^{-1}$.

\input{Tables/sart2c}

%% file: Tables/sart2a.tex
\begin{table}[H]
\centering
\caption{True error $e^1_\text{avg}(\tf)$, asymptotic error estimate $e_\text{avg}^\text{1 asymp}$, and strict error bound, $e_\text{avg}^\text{1 UB}$,  for the first-order lumped approximation of $u_\avg(t)$, $\uavgL{1}$,  as a function of $B$ (and $\Bdunk$) for domain SART-2. The time interval is given by for $[0,\tf]$ for $\tf = 2\taueq$; note that, from the definition \eqref{eq:eavgLdef}, $e^1_\text{avg}(\tf)$ is not the true error at time $\tf$, but rather the maximum of the true error over the time interval $[0,\tf]$.}
\label{tab:NumRes2a}
\begin{tabular}{l|c|l|c|c}
\multicolumn{1}{c|}{$B$} & $\Bdunk\equiv B/\gamma$ & $e^1_\text{avg}(\tf)$ & $e_\text{avg}^\text{1 asymp}$ & $e_\text{avg}^\text{1 UB}$ \\
\hline\hline
$\sci{1}{-3}$ & $\sci{1.51}{-5}$ & $\sci{8.89}{-4}$ & $\sci{8.97}{-4}$ & $\sci{2.47}{-2}$\\
$\sci{2}{-3}$ & $\sci{3.03}{-5}$ & $\sci{1.76}{-3}$ & $\sci{1.79}{-3}$ & $\sci{3.49}{-2}$\\
$\sci{5}{-3}$ & $\sci{7.57}{-5}$ & $\sci{4.27}{-3}$ & $\sci{4.49}{-3}$ & $\sci{5.52}{-2}$\\
$\sci{1}{-2}$ & $\sci{1.51}{-4}$ & $\sci{8.14}{-3}$ & $\sci{8.97}{-3}$ & $\sci{7.81}{-2}$\\
$\sci{2}{-2}$ & $\sci{3.03}{-4}$ & $\sci{1.49}{-2}$ & $\sci{1.79}{-2}$ & $\sci{1.10}{-1}$\\
$\sci{5}{-2}$ & $\sci{7.57}{-4}$ & $\sci{2.93}{-2}$ & $\sci{4.49}{-2}$ & $\sci{1.75}{-1}$\\
$\sci{1}{-1}$ & $\sci{1.51}{-3}$ & $\sci{4.31}{-2}$ & $\sci{8.97}{-2}$ & $\sci{2.47}{-1}$\\
$\sci{2}{-1}$ & $\sci{3.03}{-3}$ & $\sci{5.59}{-2}$ & $\sci{1.79}{-1}$ & $\sci{3.49}{-1}$\\
$\sci{5}{-1}$ & $\sci{7.57}{-3}$ & $\sci{6.83}{-2}$ & $\sci{4.49}{-1}$ & $\sci{5.52}{-1}$\\
$\sci{1}{0}$  & $\sci{1.51}{-2}$ & $\sci{7.46}{-2}$ & $\sci{8.97}{-1}$ & $\sci{7.81}{-1}$\\
\end{tabular}
\end{table}

%% file: Tables/sart2b.tex
\begin{table}[H]
	\centering
	\caption{True error $e^{2\text{P}}_\text{avg}(\tf)$ and asymptotic error estimate $e_\text{avg}^\text{2P asymp}$ for the second-order lumped approximation of $u_\avg(t)$, $\uavgL{2}$,  as a function of $B$ (and $\Bdunk$) for domain SART-2. The time interval is given by $[0,\tf]$ for $\tf = 2\taueq$; note that, from the definition \eqref{eq:eavgLdef}, $e^{2\text{P}}_\text{avg}(\tf)$ is not the true error at time $\tf$, but rather the maximum of the true error over the time interval $[0,\tf]$.}\label{tab:EavgL2_sart2}
	\begin{tabular}{ l | c | c | l }
	\multicolumn{1}{c|}{$B$} & $\Bdunk\equiv B/\gamma$ & $\eavgLP(\tf)$ & \multicolumn{1}{|c|}{$\eavgL{2\text{P\,asymp}}$} \\[.1em]
	\hline
	\hline
$\sci{1}{-3}$ & $\sci{1.51}{-5}$ & $\sci{9.18}{-6}$ & $\sci{1.38}{-5}$ \\
$\sci{2}{-3}$ & $\sci{3.03}{-5}$ & $\sci{3.72}{-5}$ & $\sci{5.52}{-5}$ \\
$\sci{5}{-3}$ & $\sci{7.57}{-5}$ & $\sci{2.26}{-4}$ & $\sci{3.45}{-4}$ \\
$\sci{1}{-2}$ & $\sci{1.51}{-4}$ & $\sci{8.56}{-4}$ & $\sci{1.38}{-3}$ \\
$\sci{2}{-2}$ & $\sci{3.03}{-4}$ & $\sci{3.07}{-3}$ & $\sci{5.52}{-3}$ \\
$\sci{5}{-2}$ & $\sci{7.57}{-4}$ & $\sci{1.43}{-2}$ & $\sci{3.45}{-2}$ \\
$\sci{1}{-1}$ & $\sci{1.51}{-3}$ & $\sci{3.90}{-2}$ & $\sci{1.38}{-1}$ \\
$\sci{2}{-1}$ & $\sci{3.03}{-3}$ & $\sci{9.09}{-2}$ & $\sci{5.52}{-1}$ \\
$\sci{5}{-1}$ & $\sci{7.57}{-3}$ & $\sci{2.18}{-1}$ & $\sci{3.45}{ 0}$ \\
$\sci{1}{ 0}$ & $\sci{1.51}{-2}$ & $\sci{3.54}{-1}$ & $\sci{1.38}{ 1}$ 
	\end{tabular}
\end{table}

%% file: Tables/sart2c.tex
\begin{table}[H]
	\centering
		\caption{ \ Relative error and associated asymptotic relative error estimator for the second-order approximation of $u_\Delta$, $\uDeltaL$, as a function of $B$ (and $\Bdunk$) for domain SART-2. The time interval is given by  $[0,\tf]$ for $\tf = 2\taueq$, and we choose $t_0 = 0.2 \taueq$. The third column is the true relative error, $\eDeltaLr(t_0,\tf;B)$. The fourth and fifth columns derive from the asymptotic relative error estimator provided in \eqref{eq:eDeltarelbound}: the fourth column is the long-time contribution to the asymptotic relative error estimator, and the fifth column is the full asymptotic relative error estimator --- the sum of the transient and long-time contributions.}\label{tab:EDeltaL_sart2}
         \begin{tabular}{ l | c | c | l | l}
	\multicolumn{1}{c|}{$B$} & $\Bdunk\equiv B/\gamma$ & $\eDeltaLr(2;0.1,B)$ & \multicolumn{1}{|c|}{$C_1 \Bdunk$} & \multicolumn{1}{|c}{$(C_0/(B\gamma t_0) + C_1)\Bdunk$} \\[.1em]
	\hline
	\hline
   $\sci{1}{-3}$ & $\sci{5.48}{-5}$ & $\sci{1.36}{-3}$ & $\sci{1.81}{-3}$ & $\quad\quad\sci{3.52}{-3}$ \\
   $\sci{2}{-3}$ & $\sci{1.10}{-4}$ & $\sci{2.71}{-3}$ & $\sci{3.61}{-3}$ & $\quad\quad\sci{7.03}{-3}$ \\
   $\sci{5}{-3}$ & $\sci{2.74}{-4}$ & $\sci{6.74}{-3}$ & $\sci{9.03}{-3}$ & $\quad\quad\sci{1.76}{-2}$ \\
   $\sci{1}{-2}$ & $\sci{5.48}{-4}$ & $\sci{1.33}{-2}$ & $\sci{1.81}{-2}$ & $\quad\quad\sci{3.52}{-2}$ \\
   $\sci{2}{-2}$ & $\sci{1.10}{-3}$ & $\sci{2.64}{-2}$ & $\sci{3.61}{-2}$ & $\quad\quad\sci{7.03}{-2}$ \\
   $\sci{5}{-2}$ & $\sci{2.74}{-3}$ & $\sci{8.32}{-2}$ & $\sci{9.03}{-2}$ & $\quad\quad\sci{1.76}{-1}$ \\
   $\sci{1}{-1}$ & $\sci{5.48}{-3}$ & $\sci{2.08}{-1}$ & $\sci{1.81}{-1}$ & $\quad\quad\sci{3.52}{-1}$ \\
   $\sci{2}{-1}$ & $\sci{1.10}{-2}$ & $\sci{3.80}{-1}$ & $\sci{3.61}{-1}$ & $\quad\quad\sci{7.03}{-1}$ \\
   $\sci{5}{-1}$ & $\sci{2.74}{-2}$ & $\sci{5.82}{-1}$ & $\sci{9.03}{-1}$ & $\quad\quad\sci{1.76}{ 0}$ \\
   $\sci{1}{ 0}$ & $\sci{5.48}{-2}$ & $\sci{6.79}{-1}$ & $\sci{1.81}{ 0}$ & $\quad\quad\sci{3.52}{ 0}$
	\end{tabular}
\end{table}

%% file: lump.tex
\section{Lumped Formulation: Physical Interpretation}\label{sec:lump}

In this section we apply engineering thermal resistance concepts to (i) understand the physical significance of $\phi$, and (ii) develop our small-Biot approximations. Note in this section there is no pretense of rigor: rigorous error estimates for the lumped approximations --- which serve as justifications for our choices --- are provided in Section \ref{sec:result}. Here in Section \ref{sec:lump} we shall adorn variables with $\hat\cdot$ to emphasize the formal nature of the discussion.  The main result of this section is the engineering motivation for our three lumped approximations: (i) the (first-order) classical lumped approximation, $\uavgL{1}$, \eqref{eq:u1avg}; (ii) the second-order Pad\'e lumped approximation, $\uavgLP$ \eqref{eq:u2avg}, and (iii) a second-order domain-boundary average Pad\'e approximation, $\uDeltaL$ \eqref{eq:u2Delta}.

The critical ingredient is the body average resistance $\dv{\calR}_\avg$: the effective resistance between the body average temperature and the surface temperature. We will represent this body average resistance through an effective conduction length scale $\Lcond$ based on the standard intrinsic length scale $\calL$ but now corrected for second-order effects. We will then move to a thermal circuit analogy to complete the derivations. As we shall see, $\phi$ will play a critical role. 

\subsection{Body Average Resistance} \label{sec:resistance}

We first introduce a conduction problem in which (non-uniform) heat generation balances (uniform) heat flux ($q$) at the boundary:
\begin{align}
-\dv{\bnabla}\cdot(\dv{k}\,\dv{\bnabla} \,\dv{T}^*)\ &=\ \frac{\dv{\rho c}}{\dashint_{\dv{\Omega}} \dv{\rho c}}\,\dv{q}\,\dv{ \gamma} \text{ in } \dv{\Omega}  \label{eq:strong1}, \\
\dv{k} \partial_\dv{n} \dv{T}^*\ &=-\dv{q} \text{ on } \dv{\partial\Omega} \label{eq:strong2}, \\
\dashint_{\dv{\Omega}} \frac{\dv{\rho c}}{\dashint_{\dv{\Omega}} \dv{\rho c}}(\dv{T}^* - \dv{T}_{\text{ref}})\ &= 0 \label{eq:strong3}\,,
\end{align}
where $\dv{T}_\text{ref}$ is a reference temperature. We can demonstrate that equations \eqref{eq:strong1} -- \eqref{eq:strong3} are solvable: we integrate \eqref{eq:strong1} over $\dv{\Omega}$, apply \eqref{eq:strong2}, and recall the definition of $\dv{\gamma}$ \eqref{eq:gammadef} of Section \ref{sec:form}. We now define the associated body average thermal resistance as
\begin{align}
\dv{\mathcal{R}}_{\text{avg}} &\ \equiv\ \dfrac{ \dashint_\dv{\Omega} \frac{\dv{\rho c}}{\dashint_{\dv{\Omega}} \dv{\rho c}}\,\dv{T}^* - \dashint_{\dv{\partial\Omega}} \dv{T}^*}  { \dv{q} |\dv{\partial\Omega} | } \\
&\ =\ \dfrac{ \dashint_\dv{\Omega} \frac{\dv{\rho c}}{\dashint_{\dv{\Omega}} \dv{\rho c}}\, (\dv{T}^* - \dv{T}_{\text{ref}}) - \dashint_{\dv{\partial\Omega}} (\dv{T}^* - \dv{T}_{\text{ref}}) } { \dv{q} |\dv{\partial\Omega }| } \label{eq:seceq},
\end{align}
in which the second equality \eqref{eq:seceq} follows since $\dv{T}_\text{ref}$ is a constant.

We now introduce the nondimensional temperature
\[ \hat{u} \equiv \dfrac{ \dv{T}^* - \dv{T}_{\text{ref}} }{ (\dv{\ell} \dv{q})/\dv{k}_{\inf}} \;,  \] 
substitution of which into equations \eqref{eq:strong1} -- \eqref{eq:strong3} yields
\begin{align}
-\bnabla\cdot(\kappa\bnabla \hat{u}) &= \sigma \gamma \text{ in } \Omega \label{eq:weakeq}, \\
\kappa \partial_n \hat{u} &= -1 \text{ on } \partial\Omega, \\
\int_\Omega \sigma \hat{u} &= 0 \label{eq:zeromean}.
\end{align}
The average resistance then takes the form
\begin{align}
\dv{\mathcal{R}}_{\text{avg}} &\equiv \dfrac{ \dv{\ell }} {\dv{k}_{\inf} |\dv{\partial\Omega} |} \left(-\dashint_{\partial\Omega} \hat{u}\right), \label{eq:ravdef2m1}
\end{align}
in which we have taken advantage of \eqref{eq:zeromean} to eliminate the integral over $\dv{\Omega}$. We now choose to write
\begin{align}
\dv{\mathcal{R}}_{\text{avg}} = \dfrac{\dv{\mathcal{L}}_{\text{cond}}}{\dv{k}_{\inf} |\dv{\partial\Omega}|}\; \label{eq:ravgdef2}
\end{align}
which from \eqref{eq:ravdef2m1} then yields 
\begin{align}
\mathcal{L}_\text{cond}\equiv\dv{\mathcal{L}}_{\text{cond}}/\dv{\ell} &= -\dashint_{\partial\Omega} \hat{u} \label{eq:lscale}
\end{align}
as our implicit definition of the nondimensional conduction resistance length scale. Recall that $\hat{u}$ satisfies \eqref{eq:weakeq}-\eqref{eq:zeromean}.

We now note from comparison of \eqref{eq:xi1}-\eqref{eq:xi3} and \eqref{eq:weakeq}-\eqref{eq:zeromean} that $\hat{u} = |\Omega|^{1/2} \xi$ and hence
\begin{align}
-\dashint_{\partial\Omega} \hat{u} = -|\Omega|^{1/2} \dashint_{\partial\Omega} \xi = -\dfrac{|\Omega|}{|\partial\Omega|} |\Omega|^{-1/2} |\partial\Omega| \dashint_{\partial\Omega} \xi = \text{ (from \eqref{eq:kappa_alt})  } \calL \phi  \; ,
\end{align}
or, from \eqref{eq:lscale},
\begin{align}
\mathcal{L}_{\text{cond}} &= \phi\, \calL \;.\label{eq:lcond}
\end{align}
We thus observe that $\calL_\text{cond}$, the correct conduction length scale to represent the body average thermal resistance, is the usual length scale $\calL$ amplified by $\phi$. (We recall that $\phi$ is positive, and hence $\calL_{\text{cond}}$ is a proper length.) The conduction length scale $\calL_\text{cond}$ will depend, through $\phi$, on $\sigma$, $\kappa$, and of course $\Omega$.
Our discussion here thus extends, in a quantitative fashion, the usual notion of resistance between two surfaces to the case of resistance from body ``center'' to body boundary.  
We note that there are many possible definitions of a body resistance; our particular definition ensures second-order accuracy in the resulting lumped approximations.

\subsection{Thermal Circuit}\label{sec:circuit}

\begin{figure}[H]
  \centering
  \begin{tikzpicture}[circuit ee IEC, set resistor graphic=var resistor IEC graphic]
  \tikzset{point/.style={circle,fill,inner sep=1pt}}

  \node (Tavg) at (0,0) [point,label=left:$\dv{\hat{T}}_\avg$] {};
  \node (Tpaavg) at (3,0) [point,label=above:$\dv{\hat{T}}_\paavg$] {};
  \node (Tinf) at (6,0) [point,label=right:$\dv{{T}}_\infty$] {};

  \draw (Tavg) to [resistor={info={$\dv{\mathcal{R}}_\avg$}}] (Tpaavg)
        (Tpaavg) to [resistor={info={$\dv{\mathcal{R}}_h$}}] (Tinf);

  \end{tikzpicture}
  \caption{Thermal circuit: three temperature nodes and two thermal resistances.}
  \label{fig:therm_circ}
\end{figure}
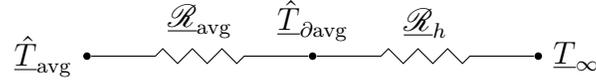

In this section, we continue to motivate the lumped approximations. To that effect, we now introduce a thermal circuit, shown in Figure \ref{fig:therm_circ}, with three nodes and associated nodal temperatures, $\dv{\hat{T}}_\avg$, $\dv{\hat{T}}_{\paavg}$, and $\dv{{T}}_{\infty}$. The nodes associated with $\dv{\hat{T}}_{\avg}$ and $\dv{\hat{T}}_{\paavg}$ are connected by a body average resistance $\dv{\mathcal{R}}_{\text{avg}}$, and the nodes $\dv{\hat{T}}_{\paavg}$ and $\dv{{T}}_{\infty}$ are connected by a standard ``heat transfer coefficient'' resistance, $\dv{\mathcal{R}}_h \equiv 1/(\dv{h}|\dv{\partial\Omega}|)$. It is important to note that $\hat\cdot$ refers to approximations to the standard quantities: $\dv{\hat{T}}_\avg$ and $\dv{\hat{T}}_\paavg$ are representative of $\dv{T}_\avg$ and $\dv{T}_\paavg$, but in general not equal to those quantities.

We may now apply standard linear circuit theory to this network. In particular, we can evaluate the heat transfer rate through our circuit as
\begin{align}
\dv{\hat{Q}} = \dfrac{\dv{\hat{T}}_\avg - \dv{{T}}_{\infty}}{\dv{\mathcal{R}}_{\text{eq}}}
\end{align}
and the temperature $\dv{\hat{T}}_{\paavg}$ as
\begin{align}
\dfrac{\dv{\hat{T}}_{\paavg} - \dv{{T}}_{\infty}}{\dv{\hat{T}}_\avg - \dv{{T}}_{\infty}} = \dfrac{\dv{\calR}_h}{\dv{\mathcal{R}}_{\text{eq}}} \label{eq:Tpaavg},
\end{align}
where $\dv{\calR}_\text{eq}$ is the (total, or) equivalent thermal resistance of the circuit,
\begin{align}
\dv{\calR}_\text{eq} \equiv \dv{\calR}_\avg + \dv{\calR}_h =\, \phi\dfrac{\dv{\calL}}{\dv{k}_{\inf} |\dv{\partial\Omega}|} + \dfrac{1}{\dv{h} |\dv{\partial\Omega}|} . \label{eq:req_sum}
\end{align}
We  observe that a large value of $\phi$ increases the body average, and hence total, resistance --- and thus decreases the heat transfer rate to the ambient; alternatively, a large value of $\phi$ decreases the fraction of the total temperature difference, $\dv{\hat{T}}_\avg - \dv{{T}}_{\infty}$, allocated to resistance $\dv{\calR}_h$ --- which thus again decreases the heat transfer rate to the ambient. Our analysis is perforce approximate, since our body average resistance assumes uniform source and boundary flux; however, we anticipate that our approximation will be quite good for small Biot number.

We now note that
\begin{align}
\Bdunkp \equiv \phi \Bdunk = \dfrac{\dv{\calR}_\avg}{\dv{\calR}_h}\; , \label{eq:Bi_dunkp}
\end{align}
which is  the usual definition of the Biot number in terms of resistance ratio but now with the improved conduction length scale provided by \eqref{eq:lcond}. We can then express $\dv{\calR}_\text{eq}$ as
\begin{align}
\dv{\calR}_\text{eq} = \dfrac{1}{\dv{h}|\dv{\partial\Omega}|}(1 + \Bdunkp) . \label{eq:R_eq}
\end{align}
Our interest in \eqref{eq:R_eq} is in the limit of small Biot number. (In the limit of large Biot number we would write instead $\dv{\calR}_\text{eq} = \phi\left[\dv{\calL}/(\dv{k}_{\inf} |\dv{\partial\Omega}|)\right] \left[1 + (\Bdunkp)^{-1}\right]$.)

We now turn to the time-dependent problem. We associate to the $\dv{\hat{T}}_\avg$ node a heat capacitance $(\dashint_{\dv{\Omega}}\dv{{\rho c}})|\dv{\Omega}|$ which then directly yields
\begin{align}
\left(\dashint_{\dv{\Omega}}\dv{{\rho c}}\right)|\dv{\Omega}| \partial_\dv{t} \dv{\hat{T}}_\avg = -\dfrac{\dv{\hat{T}}_\avg - \dv{{T}}_{\infty}}{\dv{\calR}_\text{eq}} \; .\label{eq:dim_lump_ivp}
\end{align}
We now define $\hat{u}_\avg \equiv (\dv{\hat{T}}_\avg - \dv{{T}}_{\infty})/(\dv{T}_\text{i} - \dv{{T}}_{\infty})$ and nondimensionalize \eqref{eq:dim_lump_ivp} to obtain
\begin{align}
\partial_t \hat{u}_\avg = -\dfrac{B \gamma}{1 + \phi B/\gamma}  \hat{u}_\avg \;, \quad t > 0 \label{eq:lump_ivp}
\end{align}
subject to $\hat{u}_\avg(t = 0) = 1$ from our uniform initial condition. 

\subsection{Small-Biot Approximations}\label{sec:approx_phys}

The solutions to \eqref{eq:lump_ivp} can now produce our small-$B$ approximations. First, to obtain the first-order lumped approximation for the domain average QoI, we simply set $\phi = 0$: this directly yields an exponential in time, $\uavgL{1}$, with time constant $\hat{\tau}^1 = 1/(B\gamma)$ ($= \taueq$).  Second, to obtain the second-order lumped approximation for the domain average QoI, we now retain $\phi$: this directly yields an exponential in time, $\uavgL{2}$, with time constant $\hat{\tau}^2 = (1 + \phi B/\gamma)/(B\gamma)$; the Pad\'e approximation naturally arises from our simple thermal circuit. We note that $\hat{\tau}^2 > \hat{\tau}^1$ (since $\phi > 0$), consistent with our earlier observation that the heat transfer rate, $\dv{\hat{Q}}$, decreases with increasing $\phi$. Finally, to obtain the second-order lumped approximation for the domain-boundary QoI, we now appeal to \eqref{eq:Tpaavg}: we directly obtain $\uDeltaL \equiv \frac{\phi B}{\gamma}/(1 + \frac{\phi B}{\gamma})$ (independent of time); again, the Pad\'e approximation naturally arises from our simple thermal circuit. Note that $\uDeltaL$ will only be accurate for small $B$, however, and thanks to the Pad\'e form, $0 \le \uDeltaL \le 1$, and furthermore --- as should be the case --- $\uDeltaL \to 1$ as $B \to \infty$.

%% file: phi.tex
\section{Stability and Classification for $\phi$}\label{sec:xi}

In this section we present two frameworks, with accompanying computational examples, for economization of the evaluation of $\phi$. Our emphasis is on understanding and predicting the dependence of $\phi$ on the spatial domain, $\Omega$. The first framework, empirical, is related to stability: we present a distance between domains with respect to which we claim that $\phi$ is Lipschitz continuous; we may then apply the distance to (say) our dictionary of canonical shapes, Table \ref{table:canonicalchi}, to ``nearby'' geometries. The second framework is related to classification: we present a sufficient (but not necessary) condition on the spatial domain for large $\phi$ --- hence identifying cases in which the textbook error criterion is not valid. 

We note that all of our examples in this section are in two space dimensions, $\Omega \subset \RR^2$; however, we know from property P4, tensorization, that our results for a two-dimensional domain can be directly extended to associated extruded three-dimensional domains. Finally, in this section we consider only uniform properties, $\sigma = \kappa = 1$, and hence $\phi$ will depend only on (the shape of) $\Omega$. Recall that, thanks to property P3, we can extend our results to variable $\kappa$ if we accept an upper bound for, rather than an exact value of, $\phi$ --- as might be appropriate, for example, in the error bound for our first-order lumped approximation.

As a secondary objective of this work, we have developed an efficient finite element procedure for the calculation of $\phi$. We first pose the constrained elliptic problem \eqref{eq:xi1} -- \eqref{eq:xi3} as a saddle problem; the latter (largely) preserves sparsity. We then invoke an adaptive $\mathbb{P}_2$ finite element procedure informed by {\it a posteriori} error estimators for $\xi$ (and subsequently $\phi$). The discrete equations are solved by the standard Matlab sparse LU procedure.  The finite element approximation converges very rapidly. And the error estimators ensure that the finite element error does not compromise any of our conclusions: in all cases, the finite element error in $\phi$ is much (much) less than $\phi$.

\subsection{Stability}

We now introduce a distance: given two domains in $\RR^2$, $\Omega_1$ and $\Omega_2$, we define
\begin{align}
\DD(\Omega_1,\Omega_2) \equiv \dfrac{ C_1 d_{\mathrm{Hausdorff}}(\Omega_1,\Omega_2) + C_2 \left| |\partial\Omega_1| - |\partial\Omega_2| \right|}{\max(\pcalD(\Omega_1),\pcalD(\Omega_2))}\; ,\label{eq:distdef}
\end{align}
where $d_{\mathrm{Hausdorff}}$ is the Hausdorff distance \cite{RockWets98}, $C_1$ and $C_2$ are positive real constants such that $C_1+C_2=1$, and $\pcalD(\Omega)$ is the diameter of $\Omega$. We also recall that $|\partial\Omega|$ refers to the measure of the boundary of $\Omega$; in the current two-dimensional context, $|\partial\Omega|$ is the perimeter. We note that, in two space dimensions, $\DD$ is invariant with respect to spatial scale. We should also emphasize that, since $\phi$ is invariant to translation, rotation, and dilation, we may translate, rotate, or dilate (say) domain $\Omega_2$ with respect to domain $\Omega_1$, in principle to minimize $\DD$; only shape matters.

We now claim that $\phi$ is Lipschitz continuous with respect to $\DD$: $|\phi_1- \phi_2|\leq K_L \DD(\Omega_1,\Omega_2)$ for some positive constant $K_L$, where $\phi_1$ and $\phi_2$ are the values of $\phi$ associated with respectively $\Omega_1$ and $\Omega_2$. Two disclaimers: we have no proof of this conjecture, apart from computational (empirical) evidence; we can not provide the continuity constant $K_L$. We now present some (representative) computational justification, which also illustrates the practical application of $\DD$.

\begin{figure}[H]
  \centering
  \subfloat[]{
  \centering
  \includegraphics[width=0.48\textwidth]{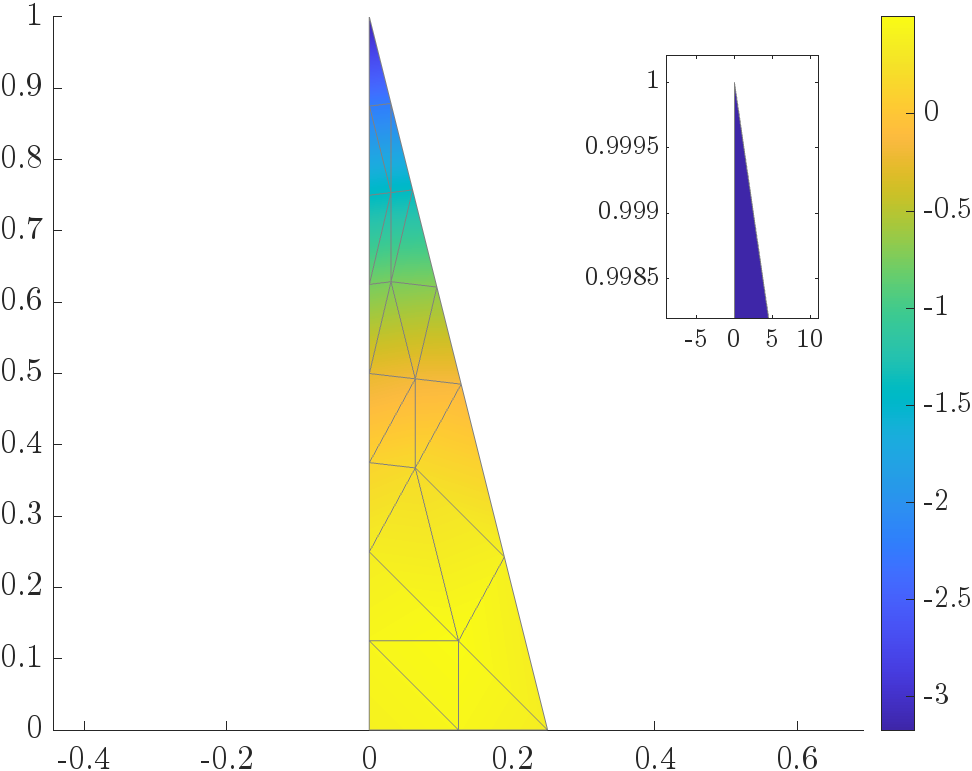}
  }
  \hfill
  \subfloat[]{
  \centering
  \includegraphics[width=0.48\textwidth]{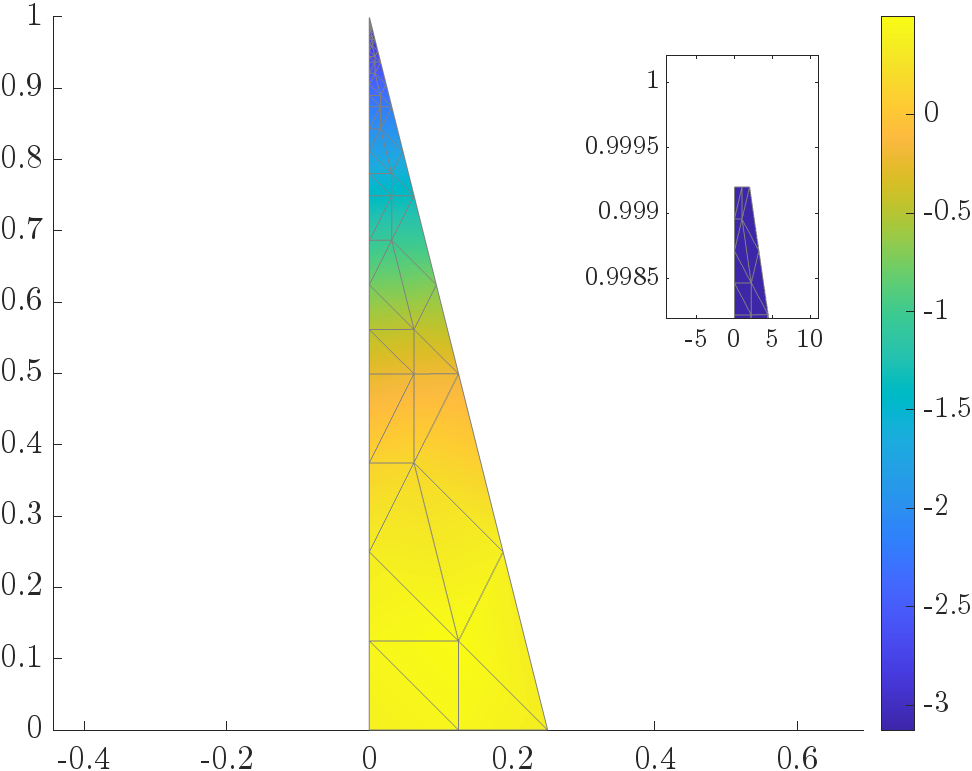}
  }\\
  \subfloat[]{
  \centering
  \includegraphics[width=0.48\textwidth]{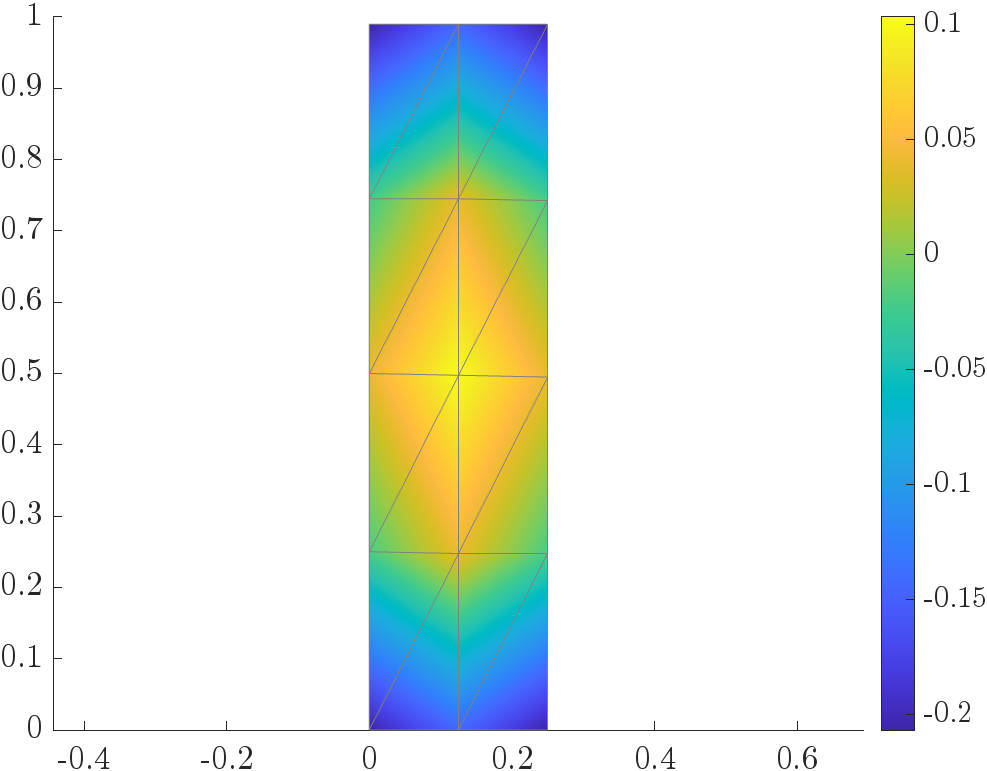}
  }
  \hfill
  \subfloat[]{
  \centering
  \includegraphics[width=0.48\textwidth]{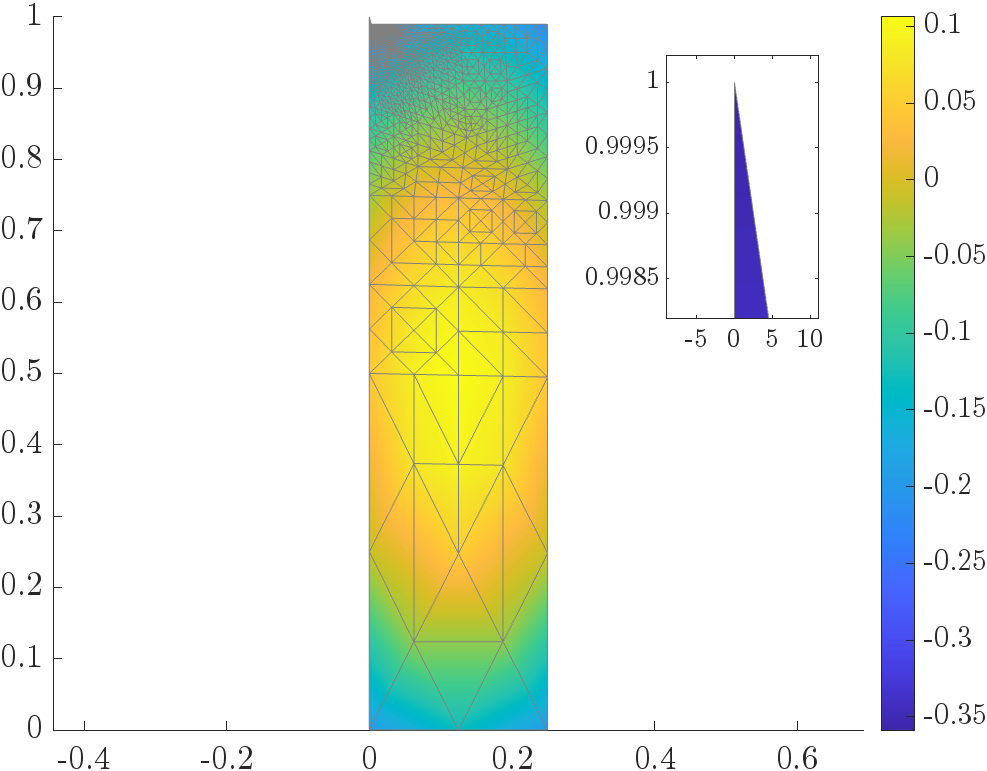}
  }
  \caption{Spatial domains for illustration of stability: a) SART-1, b) SARTC, c) RECT, and d) RECTSART. In all cases we present the geometry, the final adapted finite element mesh, and the field $\xi$ (on the basis of which $\phi$ is evaluated).} \label{fig:3phis}
\end{figure}

We consider four geometries: SART-1, SARTC, RECT, and RECTSART, presented respectively in Figure \ref{fig:3phis}a, \ref{fig:3phis}b, \ref{fig:3phis}c, and \ref{fig:3phis}d.  SARTC is a perturbation to an ``original'' geometry SART-1 $\equiv \Omega_{\text{right triangle}}(W = 1/4)$; for SART-1, $\phi = 9.13$. RECTSART is a perturbation to an ``original'' rectangle geometry RECT; for our rectangle (indeed, for any rectangle), $\phi = 2/3 = 0.667$. 
 SARTC, Figure  \ref{fig:3phis}b, is SART-1, Figure \ref{fig:3phis}a, with a $1\cdot 10^{-4}$ horizontal cut near the small angle (which thus creates a trapezoid). We obtain, for SARTC, $\phi = 9.06$, very close to $\phi$ for SART-1.  Note in this case the finite element error is not just small compared to $\phi$ but also small compared to the difference in $\phi$ between  SART-1 and SARTC. RECTSART,  Figure \ref{fig:3phis}d, is a rectangle, Figure \ref{fig:3phis}c, with an extremely small (acute) triangular extension on the upper left vertex. We obtain, for RECTSART, $\phi = 0.704$, quite close to $\phi$ for a rectangle. 

In all these examples, the distance $\DD$ between the original domain and the perturbed domain is small, and the difference in $\phi$ associated with the original domain and $\phi$ associated with the perturbed geometry is also small. It is also interesting to observe that the value of $\phi$ for SART-1, large compared to unity, is not due to the presence of a small angle. Finally, we note that, for the examples presented, the area difference term (second term) of $\DD$, \eqref{eq:distdef}, is not required; the Hausdorff term (first term) suffices. However, there are cases in which we do require the area difference term: in \cite{KandLBandP} [Section 5.4] we present the example of  a gear-like domain with increasingly oscillatory boundary; the gear domain approaches a disk in Hausdorff distance, but $\phi$ for the gear diverges --- and hence does not approach $\phi = 1/2$ for the disk.

We can draw a practical conclusion from our empirical evidence and associated observations: for a domain which, with suitable registration, is reasonably close to a canonical domain (slab, cylinder, or sphere), we anticipate that $\phi$ should be $\pcalO(1)$ --- and hence that the textbook error criterion, $\Bdunk$ small, should suffice.

\subsection{Classification: A Sufficient Feature}

We next introduce a geometric feature from which we can form a sufficient condition for $\phi$ large. We consider here homogeneous properties. For two-dimensional domains, $\Omega^{[2]}$, the feature is defined as
\begin{align}
\FF(\Omega^{[2]}) \equiv  \frac{\pi|\partial\Omega^{[2]}|^2}{8|\Omega^{[2]}|^3}\left(\text{InRadius}(\Omega^{[2]})\right)^4\,  ,\label{eq:feat2def}
\end{align}
where \text{InRadius}$(\Omega^{[2]})$ is the radius of the largest disk which can be contained within $\Omega^{[2]}$.
For three-dimensional domains, $\Omega^{[3]}$, this feature is now defined as
\begin{align}
\FF(\Omega^{[3]}) \equiv  \frac{4\pi|\partial\Omega^{[3]}|^2}{45|\Omega^{[3]}|^3}\left(\text{InRadius}(\Omega^{[3]})\right)^5\,  .\label{eq:feat2def}
\end{align}
Note that $\FF$ is scale invariant, and hence depends only on shape. This feature is actually a lower bound for $\phi$: $\phi \ge \FF$ \cite{KandLBandP} [Proposition 5.8]. Therefore, our sufficient condition for $\phi$ large is $\FF$ large. This condition can serve, in practice, as a simple discriminator: if the condition is satisfied, then the small-$B$ approximation will incur large error unless $\Bdunk$ is indeed {\em very} small. 

\begin{figure}[H]
  \centering
  \includegraphics[width=\textwidth,trim={0 11cm 0 11cm},clip]{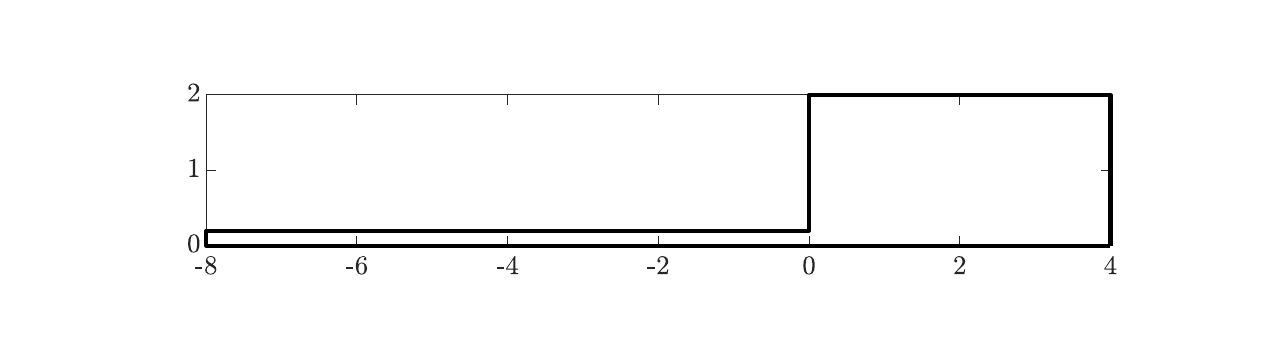}
  \caption{A ``finned-block'' domain $\Omega \subset \RR^2$ for illustration of our feature $\FF$ and associated sufficient condition. The domain comprises two subdomains: a $4\times 2$ rectangle (the block, to the right) and an $L\times H$ rectangle (the fin, to the left).} \label{fig:fin_eg}
\end{figure}

We have examined our sufficient condition over a large set of parameterized domains. We consider here
the finned-block domain shown in Figure \ref{fig:fin_eg}. The domain comprises two subdomains: a 4 $\times$ 2 rectangle (the block, to the right) and an $L \times H$ rectangle (the fin, to the left). It is readily shown that, as $L\to \infty$ while $LH$ bounded, $\FF \rightarrow \infty$, in which limit our sufficient condition is clearly satisfied. We consider the particular case $L = 8.0$ and $H = 0.2$, for which $\FF = 5.56$ and $\phi = 216$, with FE error estimate of $2.43\cdot 10^{-3}$ . In fact, $L$ (or $\FF$) need not be too large for the value of $\phi$ to be large compared to unity;  our example thus confirms our sufficient condition.

To conclude, we interpret the finned-block result physically. There are two perspectives. In the first perspective, we note that if $\FF \gg 1$, then the actual conduction length scale relevant to the internal thermal resistance of the body, clearly \text{InRadius}$(\Omega)$, may be much larger than  $\calL \equiv |\Omega|/|\partial\Omega|$. The latter is the length scale which appears in $\taueq$, \eqref{eq:taueqdef}, and thus our time constant $\taueq$  will be too small; to correct, $\phi$ inflates $\calL$ to $\calL_{\text{cond}}$, \eqref{eq:lcond}, thereby increasing the resistance and hence also the time constant. In the second perspective, we note that the temperature in the thin fin part of the domain will very quickly equilibrate, which will thus reduce the temperature difference between the fin boundary and $T_\infty$, which will in turn reduce the heat transfer rate and increase the time constant with respect to the lumped prediction --- based on uniform boundary temperature; this local reduction in the  boundary temperature is implicitly reflected in $\phi$ (through the increased length scale and resistance) as regards both $e^{\text{1 asymp}}_\avg$ and also $\uDeltaL \approx u_\Delta$.